\documentclass[a4paper,10pt]{article}
\usepackage{amsmath,amsbsy,amsfonts,amssymb,amsthm,
amscd,amsxtra,amstext,amsopn,dsfont,yfonts}
\usepackage[english]{babel}
\usepackage{psfrag}
\usepackage{indentfirst,fancyhdr,float}
\usepackage{graphicx,epsfig}
\usepackage{caption}
\usepackage{latexsym}
\usepackage{amsfonts}
\usepackage{amsmath}
\usepackage{mathrsfs}
\usepackage[all]{xy}

\def\bes{\begin{equation*}}
\def\ees{\end{equation*}}
\def\ba{\begin{aligned}}
\def\ea{\end{aligned}}
\def\be{\begin{equation}}
\def\ee{\end{equation}}
\def\bc{\begin{cases}}
\def\ec{\end{cases}}

\newcommand{\zerarcounters}{\setcounter{equation}{0}}

\theoremstyle{plain}

\theoremstyle{definition}

\let\al=\alpha \let\b=\beta    \let\g=\gamma  \let\d=\delta     
\let\e=\varepsilon
        \let\ka=\kappa     
\let\l=\lambda
\let\m=\mu    \let\n=\nu           \let\p=\pi        

\let\s=\sigma \let\t=\tau         
\let\chi=\chi
   \let\o=\omega

\def\PP{{\cal P}}\def\MM{{\cal M}} \def\VV{{\cal V}}
\def\FF{{\cal F}}\def\calH{{\cal H}}
\def\TT{{\cal T}}\def\calB{{\cal B}}
\def\RR{{\cal R}}
\def\DD{{\cal D}}\def\calA{{\cal A}}\def\GG{{\cal G}}\def\SS{{\cal S}}

\def\gotE{{\mathfrak E}}\def\FFF{{\mathfrak F}}
\def\GGG{{\mathfrak G}}

\def\gotN{{\mathfrak N}}

\def\gotR{{\mathfrak R}}

\newcommand{\ooo}{\boldsymbol{\omega}}
\newcommand{\aaa}{\boldsymbol{\alpha}}
\newcommand{\bbb}{\boldsymbol{\beta}}
\newcommand{\nnn}{\boldsymbol{\nu}}
\newcommand{\mmm}{\boldsymbol{\mu}}

\newcommand{\qqq}{\boldsymbol{q}}
\newcommand{\ppp}{\boldsymbol{p}}
\newcommand{\ccc}{\boldsymbol{c}}
\newcommand{\pps}{\psi}
\newcommand{\AAA}{\boldsymbol{A}}
\newcommand{\BBB}{\boldsymbol{B}}

\newcommand{\zerot}{0}
\newcommand{\zerooo}{\boldsymbol{0}}
\newcommand{\zerov}{\underline 0}

\def\Val{{\rm Val}}

\def\tilde#1{\widetilde{#1}}
\def\hat#1{\widehat{#1}}

\def\omt{{\omega}}
\def\omv{{\underline\omega}}
\def\nut{{\nu}}
\def\nuv{{\underline\nu}}

\def\to{\rightarrow}

\def\io{{\infty}}

\def\rme{{\rm e}}
\def\rmi{{\rm i}}

\def\RRR{{\mathbb R}}

\def\ZZZ{{\mathbb Z}} 

\def\TTT{{\mathbb T}}
\def\NNN{{\mathbb N}}
\def\id{{\mathds 1}}

\def\qed{\hfill\raise1pt\hbox{\vrule height5pt width5pt depth0pt}}

\def\ins#1#2#3{\vbox to0pt{\kern-#2 \hbox{\kern#1 #3}\vss}\nointerlineskip}

\begin{document}

\title{\bf Quasi-periodic motions in dynamical systems.\\
Review of a renormalisation group approach}

\author{
\bf Guido Gentile
\vspace{2mm}
\\ \small 
Dipartimento di Matematica, Universit\`a di Roma Tre, Roma,
I-00146, Italy.
\\ \small
E-mail: gentile@mat.uniroma3.it
}

\date{}

\maketitle

\begin{abstract}
Power series expansions naturally arise whenever solutions
of ordinary differential equations are studied in the regime
of perturbation theory. In the case of quasi-periodic solutions
the issue of convergence of the series is plagued of the so-called
small divisor problem. In this paper we review a method
recently introduced to deal with such a problem, based
on renormalisation group ideas and multiscale techniques.
Applications to both quasi-integrable Hamiltonian systems
(KAM theory) and non-Hamiltonian dissipative systems are discussed.
The method is also suited to situations in which the perturbation
series diverges and a resummation procedure can be envisaged, leading
to a solution which is not analytic in the perturbation parameter:
we consider explicitly examples of solutions which are
only $C^{\io}$ in the perturbation parameter,
or even defined on a Cantor set.
\end{abstract}

\zerarcounters
\section{Introduction}\label{sec:1}

\noindent
Consider ordinary differential equations of the form
\be
D_{\e} u = \e F(u, \omv t) ,
\label{eq:1.1} \ee
where $u=(u_{1},\ldots,u_{n})\in\RRR^{n}$,
$\e\in\RRR$ and $\omv \in \RRR^{m}$ are parameters,
called respectively the \textit{perturbation parameter}
and the \textit{frequency vector of the forcing},
$F\!:\calA \times \TTT^{m} \to \RRR^{n}$
is a real analytic function, with $\TTT=\RRR/2\p\ZZZ$
and $\calA \subset \RRR^{n}$ an open set,
and $D_{\e}$ is a differential operator possibly depending on $\e$,
\be
D_{\e} = \partial_{t}^{i_{0}} + a_{1} \e \partial_{t}^{i_{1}} +
a_{2} \e^{2} \partial_{t}^{i_{2}} + \ldots , \qquad
i_{0},i_{1},i_{2},\ldots \in \NNN .
\label{eq:1.2} \ee
In particular we shall consider explicitly the cases
$D_{\e}=\partial_{t}^{2}$ and $D_{\e}=\partial_{t}+\e\partial_{t}^{2}$.
The case $m=0$ is allowed and corresponds to $F=F(u)$.
If $m\ge 1$, we can assume without loss of generality
that the vector $\omv$ has rationally independent components.

We also assume that for $\e=0$ the \textit{unperturbed equation}
\be
D_{0} u = 0
\label{eq:1.3} \ee
admits a quasi-periodic solution $u_{0}=u_{0}(\ooo t)$,
with $\ooo\in\RRR^{p}$, $p\ge 0$,
possibly trivial (that is $p=0$, which gives a constant).
For instance, if $u$ is an angle, $u\in\TTT^{n}$,
and $D_{\e}=\partial_{t}^{2}$, one has $u_{0}=c_{0}+\Omega t$,
with $c_{0},\Omega\in\RRR^{n}$; up to a linear change of coordinates,
we can always write $\Omega=(\ooo,0,\ldots,0)$,
such that the vector $\ooo\in\RRR^{p}$, $p\le n$, has
rationally independent components.
If $u\in\RRR^{n}$ and, say, $D_{\e}=\partial_{t}+\e\partial_{t}^{2}$,
one has $u_{0}=c_{0}$, where $c_{0}\in\RRR^{n}$ is a constant vector.

We are interested in quasi-periodic solutions to (\ref{eq:1.1})
with \textit{rotation vector} (or \textit{frequency vector})
$\omt=(\ooo,\omv)\in\RRR^{d}$, $d=p+m$,
that is solutions of the form $u=u(\omt t,\e)$, with 
\be
u(\pps,\e) = \sum_{\nut\in\ZZZ^{d}}
\rme^{\rmi \nut\cdot\pps} u_{\nut}(\e) ,
\label{eq:1.4} \ee
where $\nut\cdot\pps=\n_{1}\psi_{1}+\ldots+\n_{d}\psi_{d}$ denotes
the standard scalar product. For $m=0$ we have $d=p$ and $\omt=\ooo$,
while for $p=0$ we have $d=m$ and $\omt=\omv$.
A suitable Diophantine condition will be assumed on $\omt$,
for instance the \textit{standard Diophantine condition}
\be
\left| \omt \cdot \nut \right| > \g |\nut|^{-\tau}
\quad \forall \nut\in\ZZZ^{d} \setminus\{\zerot\} ,
\label{eq:1.5} \ee
where $\g>0$ and $\t \ge d-1$ are the \textit{Diophantine constant}
and the \textit{Diophantine exponent}, respectively,
and $|\nut|=|\n_{1}|+\ldots+|\n_{d}|$. Weaker Diophantine
conditions could be considered -- see Section \ref{sec:12.1}.

The operator $D_{\e}$ acts as a multiplication operator in
Fourier space, that is
\be
(D_{\e}u)_{\nut}=\d(\omt\cdot\nut,\e)\, u_{\nut} ,
\label{eq:1.6} \ee
with $\d(x,\e)=\d_{0}(x)+
\e\d_{1}(x)+\e^{2}\d_{2}(x)+\ldots$.
For instance, if $D_{\e}=\partial_{t}^{2}$ then
$\d(\omt\cdot\nut,\e)=\d_{0}(\omt\cdot\nut)=(\rmi\omt\cdot\nut)^{2}$,
if $D_{\e}=\partial_{t}+\e\partial_{t}^{2}$ then $\d(\omt\cdot\nut,\e)=
\rmi\omt\cdot\nut \left( 1 + \rmi \e \omt\cdot\nut \right)$, and so on:
we shall explicitly take $\d(\omt\cdot\nut,\e)$ to be a
polynomial in $\e$. In other words $D_{\e}$ can
be expanded as
\be
D_{\e} = \sum_{k=0}^{k_{0}} \e^{k} D^{(k)} , \qquad D^{(0)}=D_{0} ,
\label{eq:1.7} \ee
with $k_{0}\in\NNN$ and $D^{(k)}=a_{k}\partial_{t}^{i_{k}}$, so that
$\d_{k}(\omt\cdot\nut)=a_{k}(\rmi\omt\cdot\nut)^{i_{k}}$, for $1 \le k\le k_{0}$.
The Diophantine condition (\ref{eq:1.5}) implies
\be \left| \d_{0}(\omt\cdot\nut) \right| \ge \g_{0} |\nut|^{-\tau_{0}} ,
\qquad \g_{0}=\g^{i_{0}}, \quad \tau_{0}=i_{0}\tau .
\label{eq:1.8} \ee

The problem we address here is to find a quasi-periodic solution
$u(\omt t,\e)$
to the full equation (\ref{eq:1.1}), which continues
the unperturbed solution $u_{0}(\ooo t)$, that is such that
it reduces to $u_{0}(\ooo t)$ as $\e\to0$.
This means that we look for results holding for small
values of the parameter $\e$. Hence (\ref{eq:1.1}) can
be seen as a \textit{perturbation} of the equation (\ref{eq:1.3}),
and this explains why $\e$ is called the perturbation parameter.

More precisely we shall be interested in both the existence and
stability of such quasi-periodic solutions. In particular, our
analysis accounts for the KAM theory for quasi-integrable
Hamiltonian systems (in a special case) and the existence of
quasi-periodic attractors for strongly dissipative
quasi-periodically forced one-dimensional systems.
We can also consider discrete systems,
as opposite to the continuous ones such as (\ref{eq:1.1}).
We shall see that in both cases the existence of a
quasi-periodic solution for the dynamical system is reduced
to existence of a solution for a suitable functional equation -- see
Section \ref{sec:3.2} for more insight. Extensions to more general
systems will be briefly discussed in Section \ref{sec:12}.

The method we shall follow uses renormalisation group ideas,
and is based on techniques of multiscale analysis which are
typical of quantum field theory. The method is widely
inspired to the original work of Eliasson \cite{E2} and,
even more, to its reinterpretation given by Gallavotti \cite{Ga1}.
The deep analogy with quantum field theory was stressed and used to
full extend in subsequent papers; see for instance \cite{GM1,GM2}.

For other renormalisation group approaches
existing in the literature to the same kind of problems
considered in this review see for instance \cite{MK,BGK,Koch,KLM}.

\zerarcounters
\section{Perturbation theory and formal solutions}
\label{sec:2}

\noindent
As first attempt, we look for quasi-periodic solutions (\ref{eq:1.4})
to (\ref{eq:1.1}) in the form of \textit{formal power series}
in the perturbation parameter $\e$,
\be
u(\pps,\e) = \sum_{k=0}^{\io} \e^{k} u^{(k)}(\pps) ,
\qquad u^{(k)}(\pps) = \sum_{\nut\in\ZZZ^{d}}
\rme^{\rmi \nut\cdot\pps} u^{(k)}_{\nut} , \quad k \ge 1 ,
\label{eq:2.1} \ee
where $u^{(0)}=u_{0}$ is such that $D_{0}u_{0}=0$.
The power series expansion (\ref{eq:2.1}) will be referred to
as the \textit{perturbation series} for the quasi-periodic solution.
Perturbation series have been widely studied in the literature,
especially in connection with problems of celestial mechanics \cite{Poin}.
They are sometimes called the \textit{Lindstedt series} or
\textit{Lindstedt-Newcomb series}, from the name of the astronomers
who first studied them in a systematical way.

If we expand also $D_{\e}$ according to (\ref{eq:1.7}),
we obtain to all orders $k\ge 1$
\be
D_{0} u^{(k)} = - \sum_{p=1}^{\min\{k,k_{0}\}} D^{(p)} u^{(k-p)} +
\left[ F(u,\omv t) \right]^{(k-1)} ,
\label{eq:2.2} \ee
where $[F(u,\omv t)]^{(k)}$ means that we take the Taylor
expansion of the function $u\to F(u,\cdot)$, then we
expand $u$ in powers of $\e$ according to (\ref{eq:2.1}),
and we keep the coefficient of $\e^{k}$, that is
\be
\left[ F(u,\omv t) \right]^{(k)} =
\sum_{s=0}^{\io} \frac{1}{s!}
\partial_{u}^{s} F(u_{0}(\ooo t),\omv t)
\sum_{\substack{ k_{1},\ldots,k_{s} \ge 1 \\ 
k_{1}+\ldots+k_{s}=k}} u^{(k_{1})} \ldots u^{(k_{s})} ,
\label{eq:2.3} \ee
which for $k=0$ reads $\left[ F(u,\omv t) \right]^{(0)} =
F(u_{0}(\ooo t),\omv t)$.

If $u_{0}=\ccc_{0}+\ooo t$ (and hence $p=n$), we expand
\be
F(u_{0},\omv t) = F(\ccc_{0}+\ooo t,\omv t) =
\sum_{\nut_{0} \in \ZZZ^{d}} \rme^{\rmi\nnn_{0}\cdot \ccc_{0}}
\rme^{\rmi\nut_{0}\cdot\omt t} F_{\nut_{0}} ,
\label{eq:2.4} \ee
where $\nut_{0}=(\nnn_{0},\nuv_{0})$, so as to obtain in Fourier space
\be
\d_{0}(\omt\cdot\nut) \, u^{(k)}_{\nut} = - 
\sum_{p=1}^{\min\{k,k_{0}\}} \d_{p}(\omt\cdot\nut) \,
u^{(k-p)}_{\nut} + \left[ F(u,\omv t) \right]^{(k-1)}_{\nut} ,
\label{eq:2.5} \ee
with
\be
\left[ F(u,\omv t) \right]^{(k)}_{\nut} =
\sum_{s=0}^{\io} \frac{1}{s!}
\sum_{\substack{\nut_{0},\nut_{1},\ldots,\nut_{s} \in \ZZZ^{d} \\
\nut_{0}+\nut_{1}+\ldots+ \nut_{s}=\nut}} \left( \rmi\nnn_{0} \right)^{s}
F_{\nut_{0}} \rme^{\rmi\nnn_{0}\cdot \ccc_{0}}
\sum_{\substack{k_{1},\ldots,k_{s} \ge 1 \\ k_{1}+\ldots+ k_{s}=k}}
u^{(k_{1})}_{\nut_{1}} \ldots u^{(k_{s})}_{\nut_{s}} .
\label{eq:2.6} \ee
By the analyticity assumption on the function $F$ one has
$|F_{\nut}| \le \Xi_{0} \rme^{-\xi |\nut|}$ for suitable constants
$\xi, \Xi_{0}>0$.

If $u_{0}=c_{0}$ one has $\o=\omv$. In that case we Fourier expand
$F(u,\omv t)$ only in the argument $\omv t$,
\be
F(u,\omv t) =
\sum_{\nuv_{0} \in \ZZZ^{m}} \rme^{\rmi\nuv_{0}\cdot\omv t}
F_{\nuv_{0}}(u) ,
\label{eq:2.7} \ee
and we still obtain (\ref{eq:2.5}), but with $\nut=\nuv$ and
\be
\left[ F(u,\omv t) \right]^{(k)}_{\nuv} =
\sum_{s=0}^{\io} \frac{1}{s!}
\sum_{\substack{\nuv_{0},\nuv_{1},\ldots,\nuv_{s} \in \ZZZ^{m} \\
\nuv_{0}+\nuv_{1}+\ldots+ \nuv_{s}=\nuv}}
\partial_{u}^{s} F_{\nuv_{0}} (c_{0})
\sum_{\substack{k_{1},\ldots,k_{s} \ge 1 \\ k_{1}+\ldots+ k_{s}=k}}
u^{(k_{1})}_{\nuv_{1}} \ldots u^{(k_{s})}_{\nuv_{s}} .
\label{eq:2.8} \ee

More generally, one expands
\be
\partial_{u}^{s} F(u_{0},\omv t) = 
\sum_{\nut_{0} \in \ZZZ^{d}} \rme^{\rmi\nut_{0}\cdot \omt t}
\FF_{s,\nut_{0}} ,
\label{eq:2.9} \ee
with coefficients $\FF_{s,\nut}$ bounded as
$|\FF_{s,\nut}| \le s!\Xi_{0}\Xi_{1}^{s} \rme^{-\xi |\nut|}$
for suitable constants $\xi,\Xi_{0},\Xi_{1}>0$, so that
\be
\left[ F(u,\omv t) \right]^{(k)}_{\nut} =
\sum_{s=0}^{\io} \frac{1}{s!}
\sum_{\substack{\nut_{0},\nut_{1},\ldots,\nut_{s} \in \ZZZ^{d} \\
\nut_{0}+\nut_{1}+\ldots+ \nut_{s}=\nut}} \FF_{s,\nut_{0}} 
\sum_{\substack{k_{1},\ldots,k_{s} \ge 1 \\ k_{1}+\ldots+ k_{s}=k}}
u^{(k_{1})}_{\nut_{1}} \ldots u^{(k_{s})}_{\nut_{s}} .
\label{eq:2.10} \ee
Then, to all order $k\ge 1$ one obtains
\be
u^{(k)}_{\nut} = \frac{1}{\d_{0}(\omt\cdot\nut)} \left( - 
\sum_{p=1}^{\min\{k,k_{0}\}} \d_{p}(\omt\cdot\nut) \,
u^{(k-p)}_{\nut} + \left[ F(u,\omv t) \right]^{(k-1)}_{\nut} \right) ,
\qquad \nut\neq0 ,
\label{eq:2.11} \ee
provided the \textit{compatibility condition}
\be
0 = \left[ F(u,\omv t) \right]^{(k)}_{\zerot} ,
\label{eq:2.12} \ee
holds for all $k\ge 0$. Equations (\ref{eq:2.11}) formally
provide a recursive definition of the coefficients $u^{(k)}_{\nut}$,
$\nut\neq\zerot$, in terms of the coefficients $u^{(k')}_{\nut'}$
of lower orders $k'<k$.
Indeed, the Diophantine condition (\ref{eq:1.5}) ensures
that no denominator can be zero, and the sum over the order
and Fourier labels can be easily performed (as we shall check
explicitly later on -- see Section \ref{sec:4}).
Thus, if the compatibility condition turns out to be satisfied
to all orders for a suitable choice of the coefficients
$u^{(k)}_{\zerot}$, one has an algorithm which allows to construct
iteratively all the coefficients of the series (\ref{eq:2.1}).
In that case, we say that the equations (\ref{eq:1.1}) are
\textit{formally solvable}.

\zerarcounters
\section{Examples}
\label{sec:3}

\noindent
In this section we review some examples of physically relevant
dynamical systems which can be written in the form (\ref{eq:1.1}).
Therefore, for all such systems the problem of existence of
formal quasi-periodic solutions is reduced to that of showing
that (\ref{eq:2.11}) and (\ref{eq:2.12}) can be recursively solved.
The strategy that we shall follow in the next sections will be
to show first that a formal power series solves those equations
order by order, and, then, to study the convergence of the series
for $\e$ small enough. We shall be able to prove that either
the series converges or it can be suitably resummed so as to give
a well-defined function which solves the equation (\ref{eq:1.1}).

\subsection{Quasi-integrable Hamiltonian systems: maximal tori}
\label{sec:3.1}

\noindent
Consider the Hamiltonian
\be
\calH(\aaa,\AAA) = \frac{1}{2} \AAA^{2} + \e f(\aaa) ,
\qquad \AAA^{2} = \AAA \cdot \AAA = A_{1}^{2} + \ldots + A_{n}^{2} ,
\label{eq:3.1} \ee
where $(\aaa,\AAA) \in \TTT^{n} \times \RRR^{n}$ are \textit{angle-action}
coordinates, and $f\!:\TTT^{n}\to\RRR$ is a real analytic function.

More generally we could consider Hamiltonians of the form
\be
\calH(\aaa,\AAA) = \calH_{0}(\AAA) + \e f(\aaa,\AAA) ,
\label{eq:3.2} \ee
where $\calH_{0}$ and $f$ are both real analytic functions in
$\TTT^{n} \times \calA$, $\calA \subset \RRR^{n}$ being an open set,
with $\calH_{0}$ convex and $f$ periodic in the angles $\aaa$. 
Existence and properties of quasi-periodic solutions for
\textit{quasi-integrable systems}, that is systems
described by Hamiltonians of the form (\ref{eq:3.2}),
are the content of \textit{KAM theory} \cite{Kol,A,Mo1};
see \cite{AKN} for a review.

The advantage of taking the Hamiltonian (\ref{eq:3.1}) is that the
corresponding Hamilton equations lead to a closed equation for $\aaa$,
\be
\ddot \aaa = - \e \partial_{\aaa}f(\aaa) ,
\label{eq:3.3} \ee
which is the form (\ref{eq:1.1}) with $u=\aaa$, $m=0$,
and $D_{\e}=\partial_{t}^{2}$ (so that $\d_{0}(\omt\cdot\nut)=
(\rmi\omt\cdot\nut)^{2}$, while $\d_{k}(\omt\cdot\nut)=0$ for all $k\ge1$).

We are interested in quasi-periodic solutions of the
form (\ref{eq:1.4}) with $d=p=n$ (and $\omt=\ooo$). 
For $\e=0$ we take as unperturbed solution $u_{0}(\ooo t)=
\aaa_{0} + \ooo t$, with $\aaa_{0}\in\TTT^{d}$ arbitrary
and $\omt=\ooo \in \RRR^{d}$ satisfying
the Diophantine condition (\ref{eq:1.5}).

If we write, in agreement with (\ref{eq:2.1}),
\be
\aaa(\ooo t,\e) = \aaa_{0} + \ooo t +
\sum_{k=1}^{\io} \e^{k} \sum_{\nnn\in\ZZZ^{n}}
\rme^{\rmi\ooo\cdot\nnn t} \aaa^{(k)}_{\nnn} ,
\label{eq:3.4} \ee
we find to all orders $k\ge1$ -- see 
(\ref{eq:2.11}) and (\ref{eq:2.12}) --
\begin{subequations}
\begin{align}
\aaa^{(k)}_{\nnn} & = \frac{1}{(\omt\cdot\nnn)^{2}}
\left[ \partial_{\aaa} f(\aaa) \right]^{(k-1)}_{\nnn} ,
\qquad \nnn\neq\zerooo ,
\label{eq:3.5a} \\
\zerooo & = \left[ \partial_{\aaa} f(\aaa) \right]^{(k-1)}_{\zerooo} ,
\label{eq:3.5b}
\end{align}
\label{eq:3.5} \end{subequations}
\vskip-.3truecm
\noindent where
\be
\left[ \partial_{\aaa} f(\aaa) \right]^{(k)}_{\nnn} =
\sum_{s=0}^{\io} \frac{1}{s!}
\sum_{\substack{\nnn_{0},\nnn_{1},\ldots,\nnn_{s}\in\ZZZ^{n} \\
\nnn_{0}+\nnn_{1}+\ldots+\nnn_{s}=\nnn}}
\left( \rmi \nnn_{0} \right)^{s+1} f_{\nnn_{0}}
\sum_{\substack{ k_{1},\ldots,k_{s} \ge 1 \\ k_{1}+\ldots+k_{s}=k}}
\aaa^{(k_{1})}_{\nnn_{1}} \ldots \aaa^{(k_{s})}_{\nnn_{s}} .
\label{eq:3.6} \ee
Thus, equation (\ref{eq:3.5a}) defines recursively
the coefficients $\aaa^{(k)}_{\nnn}$ for all $k\ge1$ and
all $\nnn\neq\zerooo$, provided equation (\ref{eq:3.5b})
is satisfied for all $k\ge 1$.

The compatibility condition (\ref{eq:3.5b}) is automatically satisfied
for $k=1$, because $\left[ \partial_{\aaa} f(\aaa) \right]^{(0)}_{\nnn} =
\rmi\nnn f_{\nnn}$, which vanishes for $\nnn=\zerooo$.
It is a remarkable cancellation that the condition holds for all
$k\ge 1$ -- see Section \ref{sec:4.4} --, so implying that the
perturbation series (\ref{eq:3.4}) is well-defined to all orders.
The coefficients $\aaa^{(k)}_{\zerooo}$, $k \ge 1$, can be arbitrarily
fixed; for instance one can set $\aaa^{(k)}_{\zerooo}=\zerooo$
for all $k\ge1$ -- see Section \ref{sec:4.4}.
We shall see in Section \ref{sec:7} that for $\e$
small enough the series converges to a function analytic in
$\boldsymbol{\psi}=\ooo t$. As a consequence, 
there exists a quasi-periodic solution of the form (\ref{eq:3.4}),
analytic both in $\e$ and $\boldsymbol{\psi}$, and
parameterised by $\aaa_{0}\in\TTT^{n}$: hence such a solution
describes an $n$-dimensional invariant torus (\textit{maximal KAM torus}).

In this paper we confine ourselves to the Hamiltonian
(\ref{eq:3.1}). Note that in the case (\ref{eq:3.2})
the unperturbed solution is $u_{0}=(\aaa_{0}+\ooo(\AAA_{0}) t,\AAA_{0})$,
with $\ooo(\AAA)=\partial_{\AAA}\calH_{0}(\AAA)$,
hence it is still of the form $\ccc_{0}+\ooo t$ with $n=2p$
and $\o_{i}=0$ for $i\ge n+1$. However, strictly speaking
the Hamilton equation are not of the form (\ref{eq:1.1}),
so the analysis should be suitably adapted -- see \cite{GM2}.

We could also consider, instead of (\ref{eq:3.3}),
the more general equation
\be
\ddot{\aaa} = - \e \partial_{\aaa}f(\aaa,\omv t) ,
\label{eq:3.7} \ee
which reduces to (\ref{eq:3.3}) for $m=0$. This is still
a Hamiltonian system, with Hamiltonian
\be \calH(\aaa,\underline\al,\AAA,\underline A) =
\frac{1}{2} \AAA \cdot \AAA + \omv \cdot \underline A +
\e f(\aaa,\underline\al) .
\label{eq:3.8} \ee
The unperturbed solution to (\ref{eq:3.7}) is of the same
form as before, but in that case we look for
a quasi-periodic solution with rotation vector $\omt=(\ooo,\omv)$.

Then the angles $\underline\al$ evolve trivially as $\underline\al(t)=
\omv t$, whereas equation (\ref{eq:3.5a}) has to be replaced with
\be
\aaa^{(k)}_{\nut} = \frac{1}{(\omt\cdot\nut)^{2}}
\sum_{s=0}^{\io} \frac{1}{s!}
\sum_{\substack{\nut_{0},\nut_{1},\ldots,\nut_{s}\in\ZZZ^{d} \\
\nut_{0}+\nut_{1}+\ldots+\nut_{s}=\nut}}
\left( \rmi \nnn_{0} \right)^{s+1} f_{\nut_{0}}
\sum_{\substack{ k_{1},\ldots,k_{s} \ge 1 \\ k_{1}+\ldots+k_{s}=k-1}}
\aaa^{(k_{1})}_{\nut_{1}} \ldots \aaa^{(k_{s})}_{\nut_{s}} ,
\label{eq:3.9} \ee
while the compatibility condition for $\nut=\zerot$ reads
\be
\zerooo = \sum_{s=0}^{\io} \frac{1}{s!}
\sum_{\substack{\nut_{0},\nut_{1},\ldots,\nut_{s}\in\ZZZ^{d} \\
\nut_{0}+\nut_{1}+\ldots+\nut_{s}=\zerot}}
\left( \rmi \nnn_{0} \right)^{s+1} f_{\nut_{0}} 
\sum_{\substack{ k_{1},\ldots,k_{s} \ge 1 \\ k_{1}+\ldots+k_{s}=k-1}}
\aaa^{(k_{1})}_{\nut_{1}} \ldots \aaa^{(k_{s})}_{\nut_{s}} .
\label{eq:3.10} \ee
However, the differences with respect to the previous case
are just minor ones, as one can easily work out by himself.

\subsection{Discrete systems: the standard map}
\label{sec:3.2}

\noindent
Beside continuous dynamical systems we can consider discrete
dynamical systems (maps), such as the \textit{standard map}
\cite{Chi,LL} -- also known as Chirikov-Greene-Taylor map.

The standard map is defined by the symplectic map
from the cylinder to itself
\be
\begin{cases}
x'=x+y+\e \, \sin x , & \\ y' =y+\e\,\sin x , & \end{cases}
\label{eq:3.11} \ee
where $(x,y)\in \TTT\times \RRR$.

For $\e=0$ the motion is trivial: one has a simple rotation
$x'=x+2\pi\omega$, while $y$ is fixed to $y=2\pi\omega$,
$\omega\in\RRR$. For $\e\neq 0$ one can look for solutions
which are conjugate to a trivial rotation of some other variable
(KAM invariant curves), i.e. one can look for solutions of the form
\be
x =\al + u(\al,\e) , \qquad y=2\pi\o + v(\al,\e) ,
\label{eq:3.12} \ee
with $\al\to\al'=\al+2\pi\o$ and the
functions $u,v$ depending analytically on their arguments.
The number $\o$ will be called the \textit{rotation number}.

The functions $u,v$ are not independent from each other.
One has $v(\al,\e)=u(\al,\e)-u(\al-2\pi\o,\e)$,
as it is straightforward to check:
simply note that $x'=x+y'$ by (\ref{eq:3.11}) and express
$x'$ and $y'$ in terms of $\al$ through (\ref{eq:3.12}),
using that $\al'=\al+2\pi\o$. So one obtains a closed
equation for the \textit{conjugating function} $u$,
\be
D u(\al,\e)\equiv
u (\al+2\pi\o,\e)+u(\al-2\pi\o,\e)-2u(\al,\e) =
\e \sin (\al+u(\al,\e) ).
\label{eq:3.13} \ee
If such a solution exists and is analytic in $\e$,
then it has to be possible to expand the function $u$
as Taylor series in $\e$ and as Fourier series in $\al$.
So we are led to write, at least formally,
\be
u(\al,\e)=\sum_{\nu\in \ZZZ}\sum_{k=1}^{\infty}
\rme^{\rmi\nu\al}\e^{k} u_{\nu}^{(k)} .
\label{eq:3.14} \ee
so implying, by calling $[F(\alpha,\e)]^{(k)}_{\nu}$
the coefficient of the function $F(\alpha,\e)$
with Fourier label $\nu$ and Taylor label $k$,
according to (\ref{eq:2.3}),
\be
\d_{0}(\o\nu) u_{\nu}^{(k)}=
\left[ \sin (\al+u(\al,\e) )\right]^{(k-1)}_{\nu} ,
\qquad \d_{0}(\o\nu)=2\left[ \cos (2\pi\o\nu)-1\right] .
\label{eq:3.15} \ee
Note that, for $\o\nu$ small (mod. 1), $\d_{0}(\o\nu)\sim
\|\o\nu\|^2$, if $\|x\|=\min_{p\in\ZZZ}|x-p|$
denotes the distance of $x$ from the nearest integer.
By explicitly writing
$ \sin \al = \sum_{\nu_{0}=\pm 1} (2\rmi)^{-1} \nu_{0}\rme^{\rmi\nu_{0}\al}$
and Taylor expanding $\sin(\al+u)$ in $u$
around $u=0$, one finds, from (\ref{eq:3.13}),
\begin{subequations}
\begin{align}
u_{\nu}^{(1)} & = 
-\frac{\rmi\nu}{2\d_{0}(\o\nu)} ,
\label{eq:3.16a} \\
u_{\nu}^{(k)} & =
\frac{1}{\d_{0}(\o\nu)} \sum_{s=0}^{\io}
\sum_{\substack{ \nu_{0}+\nu_{1}+\ldots+\nu_{s}=\nu
\\ k_{1}+\ldots+k_{s}=k-1}} \frac{-(\rmi\nu_{0})^{s+1}}{s!2}
u_{\nu_{1}}^{(k_{1})}\ldots u_{\nu_{s}}^{(k_{s})} , \qquad k>1 ,
\label{eq:3.16b}
\end{align}
\label{eq:3.16} \end{subequations}
\vskip-.3truecm
\noindent for $\nu\neq0$. It is easy to check that
for all $k\ge 1$ one has $u_{\nu}^{(k)}=0$ if $|\nu|> k$,
and one can choose $u^{(k)}_{0}=0$.
Then (\ref{eq:3.16}) can be iterate
by taking into account that $k_j<k$ for any $j=1,\ldots,s$.

When working in Fourier space, the recursive equations
(\ref{eq:3.16}) look very similar to the equations (\ref{eq:3.5a}) 
for continuous systems, with $n=2$. Then, one could be tempted to
write the standard map as the stroboscopic map of
a continuous system. However, it turns out that, formally,
one should consider the singular system (known in physics as the
\textit{kicked rotator}) with Hamiltonian 
\be
\calH(\aaa,\AAA) = 2\p A_{1} + \frac{A_{2}^{2}}{2} + 2\p \e
\sum_{n\in\ZZZ} \d(\al_{1}-2\p n) \left( \cos \al_{2} - 1 \right) ,
\label{eq:3.17} \ee
where $\d$ is the delta function, and set $x=\al_{2}$
and $y=A_{2}-(\e/2)\sin\al_{2}$ \cite{GBG}. Therefore
the corresponding Hamilton equations cannot be in the
class (\ref{eq:1.1}), where smoothness was required.
Nonetheless, as we have seen, in Fourier space the
analysis is essentially the same.

\subsection{Quasi-integrable Hamiltonian systems: lower dimensional tori}
\label{sec:3.3}

\noindent
Consider the Hamiltonian
\be
\calH(\aaa,\bbb,\AAA,\BBB) =
\frac{1}{2} \AAA^{2} + \frac{1}{2} \BBB^{2} +
\e \, f(\aaa,\bbb) ,
\label{eq:3.18} \ee
where $(\aaa,\AAA)\in\TTT^{r}\times\RRR^{r}$ and
$(\bbb,\BBB)\in\TTT^{s}\times\RRR^{s}$ are
angle-action coordinates, with $r+s=n$, and
$f\!:\TTT^{n}\to\RRR$ is a real analytic function.

We can also consider the same Hamiltonian as (\ref{eq:3.1}), but
assume that $\ooo$ is a resonant vector, that is that there
exist $s$ integer vectors $\nnn_{1},\ldots,\nnn_{s}$ such that
$\ooo\cdot\nnn_{1}=\ldots=\ooo\cdot\nnn_{s}=0$. If this happens,
it is possible to perform a linear change of coordinates
such that in the new coordinates $\ooo=(\o_{1},\ldots,\o_{r},
0,\ldots,0)$, with $(\o_{1},\ldots,\o_{r})$ satisfying
a Diophantine condition in $\RRR^{r}$. Of course the
corresponding Hamiltonian would be slightly more complicated
than (\ref{eq:3.18}). For simplicity's sake
we shall confine ourselves to (\ref{eq:3.18}), 
and study the problem of existence of a quasi-periodic solution
with rotation vector $\ooo=(\o_{1},\ldots,\o_{r})$ which for $\e=0$
reduces to $u_{0}=(\aaa_{0}+\ooo t,\bbb_{0},\AAA_{0},\zerooo)$,
with $\AAA_{0}=\ooo$.

In terms of the angles $(\aaa,\bbb)$ the Hamilton equations become
\be
\begin{cases}
\ddot \aaa = - \e \partial_{\aaa} f(\aaa,\bbb) , & \\
\ddot \bbb = - \e \partial_{\bbb} f(\aaa,\bbb) , &
\end{cases}
\label{eq:3.19} \ee
so that the unperturbed solution is $(\aaa_{0}+\ooo t,\bbb_{0})$,
which is of the form $u_{0}=c_{0}+\Omega t$, with
$\Omega=(\o_{1},\ldots,\o_{r},0,\ldots,0)=(\ooo,\zerooo)$.
Hence we look for solutions of the form (\ref{eq:1.4}) with $d=p=r$.

Since $\bbb$ is expected to remain close to $\bbb_{0}$,
we Fourier expand $f$ only in the angle $\aaa$, so writing
\be
f(\aaa,\bbb)=\sum_{\nnn\in\ZZZ^{d}} \rme^{\rmi \nnn\cdot\aaa}
f_{\nnn}(\bbb) .
\label{eq:3.20} \ee
Therefore, (\ref{eq:3.19}) gives, in Fourier space,
\be
\begin{cases}
\left( \ooo\cdot\nnn \right)^{2} \aaa^{(k)}_{\nnn} =
\left[ \partial_{\aaa} f(\aaa,\bbb) \right]^{(k-1)}_{\nnn} , & \\
\left( \ooo\cdot\nnn \right)^{2} \bbb^{(k)}_{\nnn} =
\left[ \partial_{\bbb} f(\aaa,\bbb) \right]^{(k-1)}_{\nnn} , &
\end{cases}
\label{eq:3.21} \ee
with
\be
\left[ \partial_{\aaa} f(\aaa,\bbb) \right]^{(k)}_{\nnn} =
\sum_{p,q=0}^{\io} \frac{1}{p!q!} \!\!\!
\sum_{\substack{\nnn_{0},\nnn_{1},\ldots,\nnn_{p+q} \in \ZZZ^{n} \\
\nnn_{0}+\nnn_{1}+\ldots+\nnn_{p+q}=\nnn}}
\!\!\!\!\!\!\!\!\!\!\!\!\!\!\!
\left( \rmi\nnn_{0} \right)^{p+1} \partial_{\bbb}^{q}
f_{\nnn_{0}}(\bbb_{0})
\!\!\!\!\!\!\!\!\!\!\!\!
\sum_{\substack{k_{1},\ldots,k_{p+q} \ge 1 \\ k_{1}+\ldots+k_{p+q}=k}}
\!\!\!\!\!\!\!\!\!\!\!\!
\aaa^{(k_{1})}_{\nnn_{1}} \ldots \aaa^{(k_{p})}_{\nnn_{p}}
\bbb^{(k_{p+1})}_{\nnn_{p+1}} \ldots \bbb^{(k_{p+q})}_{\nnn_{p+q}} ,
\label{eq:3.22} \ee
and an analogous expression holding for
$\left[ \partial_{\bbb} f(\aaa,\bbb) \right]^{(k)}_{\nnn}$ -- with
$\left( \rmi\nnn_{0} \right)^{p} \partial_{\bbb}^{q+1}$
instead of $\left( \rmi\nnn_{0} \right)^{p+1} \partial_{\bbb}^{q}$.

Again for $\nnn=\zerooo$ we require both
\be
\left[ \partial_{\aaa} f(\aaa,\bbb) \right]^{(k)}_{\zerooo}=\zerooo ,
\qquad 
\left[ \partial_{\bbb} f(\aaa,\bbb) \right]^{(k)}_{\zerooo}=\zerooo ,
\label{eq:3.23} \ee
for all $k\ge 0$.
We shall see that the first compatibility condition is automatically
satisfied for all values of $\aaa_{0}$ and $\aaa^{(k)}_{\zerooo}$,
$k\ge 1$, while the second one requires $\bbb_{0}$ and
$\bbb^{(k)}_{\zerooo}$, $k\ge 1$, to be suitably fixed.
This is clear already to first order, where we obtain
\be
\begin{cases}
\left( \ooo\cdot\nnn \right)^{2} \aaa^{(1)}_{\nnn} =
\rmi \nnn f_{\nnn}(\bbb_{0}) , & \\
\left( \ooo\cdot\nnn \right)^{2} \bbb^{(1)}_{\nnn} =
\partial_{\bbb} f_{\nnn}(\bbb_{0}) , &
\end{cases}
\label{eq:3.24} \ee
so that for $\nnn=\zerooo$ the first equation trivially holds,
whereas the second one fixes $\bbb_{0}$ to be such that
$\partial_{\bbb} f_{\zerooo}(\bbb_{0}) =\zerooo$, i.e. $\bbb_{0}$
must be a stationary point for the function
$f_{\zerooo}(\bbb_{0})$ (such a point always exists).

Moreover to higher orders one has
\be
\left[ \partial_{\bbb} f(\aaa,\bbb) \right]^{(k)}_{\zerooo} =
\partial_{\bbb}^{2}f_{\zerooo}(\bbb_{0}) \, \bbb^{(k)}_{\zerooo}
+ \boldsymbol{\Phi}_{k} ,
\label{eq:3.25} \ee
for a suitable function $\boldsymbol{\Phi}_{k}$
depending only on the coefficients of order strictly less
than $k$ -- see Section \ref{sec:4.4} for details.
Thus, if we further assume that the matrix
$\partial_{\bbb}^{2}f_{\zerooo}(\bbb_{0})$ be nonsingular
(\textit{nondegeneracy condition}),
then we can impose the compatibility conditions
$\left[ \partial_{\bbb} f(\aaa,\bbb) \right]^{(k)}_{\zerooo} = \zerooo$
by suitably fixing the corrections $\bbb^{(k)}_{\zerooo}$, $k \ge 1$,
to the constant part of the $\bbb$ angles.

In Section \ref{sec:4} we shall see that, at least formally,
a quasi-periodic solution parameterised by $\aaa_{0}\in\TTT^{r}$
exists for suitable values of $\bbb_{0}$. We shall see
in Section \ref{sec:8} that, even if the formal series is divergent,
however it can be suitably resummed for $\e$ small enough
so as to be given a meaning
as a well-defined function analytic in $\boldsymbol{\psi}=\ooo t$:
hence the latter describes a \textit{lower-dimensional torus}.
If the matrix $\partial_{\bbb}^{2}f_{\zerooo}(\bbb_{0})$ is
positive definite we shall say that the lower-dimensional torus
is \textit{hyperbolic} if $\e<0$ and \textit{elliptic}
if $\e>0$; in the latter case we shall see that the torus
exists only for some values of $\e$, more precisely for
$\e$ defined in a Cantor set with Lebesgue density point at the origin
(\textit{Cantorisation}) -- see Section \ref{sec:10}.

\subsection{Strongly dissipative quasi-periodically forced systems}
\label{sec:3.4}

\noindent
Consider a one-dimensional system subject to a mechanical force $g$,
in the presence of dissipation and of a quasi-periodically forcing.
The equation describing the system is the
ordinary differential equation
\be \ddot x + \g \dot x + g(x) = f(\omv t) ,
\label{eq:3.26} \ee
where $x\in\RRR$, $\g>0$ is the \textit{dissipation coefficient}
and $\omv\in\RRR^{m}$ is the frequency vector of the forcing.
We assume that both $g\!:\cal A \to \RRR$ and $f\!:\TTT^{m}\to\RRR$
are real analytic functions, with $\calA\subset\RRR$ an open set.

If the dissipation is large enough, it is natural to rewrite
(\ref{eq:3.26}) in terms of the small parameter $\e=1/\g$,
so as to obtain the equation
\be \dot x + \e \ddot x + \e g(x) = \e f(\omv t) ,
\label{eq:3.27} \ee
which is of the form (\ref{eq:1.1}) with $n=1$, $u=x$,
$D_{\e}=\partial_{t} + \e \partial_{t}^{2}$, and $F(u,\omv t)
=-g(u)+f(\omv t)$. In particular one has
$\d_{0}(\omt\cdot\nut)=\rmi\omt\cdot\nut$,
$\d_{1}(\omt\cdot\nut)=(\rmi\omt\cdot\nut)^{2}$,
$a_{1}=1$, and $\d_{k}(\omt\cdot\nut)=a_{k}=0$ for all $k\ge2$.

We look for quasi-periodic solutions $x(\omv t,\e)$ which are
analytic in $\underline{\psi}=\omv t$ and continue
the unperturbed solutions $x_{0}=c_{0}$, with $c_{0}$ constant.
Such solutions (if any) are called \textit{response solutions},
as they have the same frequency vector as the forcing.
Thus, $\omt=\omv$ and $d=m$, so that, if we write
\be
x(\omv t,\e) = c_{0} +
\sum_{k=1}^{\io} \e^{k} \sum_{\nuv\in\ZZZ^{m}}
\rme^{\rmi\nuv\cdot\omv t} x^{(k)}_{\nuv} ,
\label{eq:3.28} \ee
then we obtain
\be
\rmi\omt\cdot\nuv \, x^{(k)}_{\nuv} +
\left( \rmi\omv\cdot\nuv \right)^{2} x^{(k-1)}_{\nuv} +
[g(x)]^{(k-1)}_{\nuv} = f_{\nuv} \d_{k,1}
\label{eq:3.29} \ee
for $k\ge1$, while $x^{(0)}_{\zerov}=c_{0}$ and $x^{(0)}_{\nuv} = 0$
for $\nuv \neq \zerov$.

The first order equation gives
\be
\begin{cases}
\rmi\omv\cdot\nuv \, x^{(1)}_{\nuv} = f_{\nuv} , &
\quad \nuv\neq\zerov \\
g(c_{0}) = f_{\zerov} , \end{cases}
\label{eq:3.30} \ee
which fixes the value of $c_{0}$ (of course $f_{\zerov}$ must
belong to the range of $g$). The unperturbed solution
can be seen as a quasi-periodic solution in the
extended phase space $(x,\underline{\psi}) \in \RRR \times \TTT^{m}$,
where it looks like $(c_{0},\omv t)$, so that the full solution
is of the form $(c_{0}+X(\omv t,c_{0},\e),\omv t)$.

If $\partial_{x}g(c_{0})\neq0$, then to higher orders one has
the compatibility conditions
\be
\left[ g(x) \right]^{(k)}_{\zerov} =
\partial_{x}g(c_{0}) \, x^{(k)}_{\zerov} +
G_{k}(c_{0},x^{(1)},\ldots,x^{(k-1)}) = 0 , \qquad k \ge 1 ,
\label{eq:3.31} \ee
for a suitable function $G_{k}$ depending only on the coefficients
of orders $k'<k$. An explicit calculation gives
\be
G_{k}(c_{0},x^{(1)},\ldots,x^{(k-1)}) = 
\sum_{s=2}^{\io} \frac{1}{s!} \partial_{x}^{s} g(c_{0}) 
\sum_{\substack{\nuv_{1},\ldots,\nuv_{s} \in \ZZZ^{m} \\ 
\nuv_{1}+\ldots+\nuv_{s}=\zerov}}
\sum_{\substack{k_{1},\ldots,k_{s} \ge 1 \\ k_{1}+\ldots+k_{s}=k}}
x^{(k_{1})}_{\nuv_{1}} \ldots x^{(k_{s})}_{\nuv_{s}} .
\label{eq:3.32} \ee
Thus, (\ref{eq:3.31}) can be used to fix
the corrections $x^{(k)}_{\zerov}$, $k\ge 1$,
to the constant part of the solution $x(\omv t,\e)$.

We shall see in Section \ref{sec:9} that, under the
\textit{nondegeneracy condition} $\partial_{x}g(c_{0}) \neq0$,
the series (\ref{eq:3.28}) can be resummed for $\e$ small enough
into a function which depends analytically on $\underline{\psi}=\omv t$.

\zerarcounters
\section{Diagrammatic representation and tree formalism}
\label{sec:4}

\noindent
We have to study the recursive equations (\ref{eq:2.11})
and (\ref{eq:2.12}), with $\left[ F(u,\omv t) \right]^{(k-1)}_{\nut}$
given by (\ref{eq:2.10}). This will be done through a
diagrammatic formalism, known as the \textit{tree formalism}.

Let us assume that one can decompose
$u=(\tilde u, \hat u)$, with $\tilde u \in \RRR^{\tilde n}$ and
$\hat u \in \RRR^{\hat n}$, $\tilde n + \hat n = n$, and,
accordingly, $F=(\tilde F,\hat F)$, so that for $\nut=\zerot$ one has
\begin{subequations}
\begin{align}
[ \tilde F(u,\omv t) ]^{(k)}_{\zerot} & = \tilde 0 ,
\label{eq:4.1a} \\
[ \hat F(u,\omv t) ]^{(k)}_{\zerot} & =
A \, \hat u^{(k)}_{\zerot} + \Phi_{k} ,
\label{eq:4.1b}
\end{align}
\label{eq:4.1} \end{subequations}
\vskip-.3truecm
\noindent for a suitable nonsingular matrix $A$ and a suitable
vector $\Phi_{k}$ ($\tilde 0$ is the null vector in $\RRR^{\tilde n}$).
We shall see that in all cases considered in Section \ref{sec:3}
this holds true -- see Section \ref{sec:4.4} below.

We first introduce the trees (see also \cite{Bo,H,HP}) as
the main combinatorial and graphical objects that we shall
use in the forthcoming analysis. Then we shall provide
some rules how to associate numerical values to the trees,
so as to represent the coefficients $u^{(k)}_{\nut}$ in terms of trees.

\subsection{Trees}
\label{sec:4.1}

\noindent
A connected \textit{graph} $\GG$ is a collection of points (nodes) and
lines connecting all of them. Denote by $N(\GG)$ and $L(\GG)$ the set
of nodes and the set of lines, respectively. A path between two nodes
is the minimal subset of $L(\GG)$ connecting the two nodes. A graph is
planar if it can be drawn in a plane without graph lines crossing.

A \textit{tree} is a planar graph $\GG$ containing no closed loops.
Consider a tree $\GG$ with a single special node $v_{0}$:
this introduces a natural partial ordering on the set
of lines and nodes, and one can imagine that each line
carries an arrow pointing toward the node $v_{0}$.
We add an extra oriented line $\ell_{0}$ exiting the special
node $v_{0}$; the added line will be called the \textit{root line}
and the point it enters (which is not a node) will be called the
\textit{root} of the tree. In this way we obtain a \textit{rooted tree}
$\theta$ defined by $N(\theta)=N(\GG)$ and $L(\theta)=L(\GG)\cup\ell_{0}$.
A \textit{labelled tree} is a rooted tree $\theta$ together with a label
function defined on the sets $L(\theta)$ and $N(\theta)$.

We call \textit{equivalent} two rooted trees which can be
transformed into each other by continuously deforming the lines
in the plane in such a way that the lines do not cross each other.
We can extend the notion of equivalence also to labelled trees,
by considering equivalent two labelled trees if they can be
transformed into each other in such a way that the labels also match. 
In the following we shall deal mostly with nonequivalent labelled
trees: for simplicity, where no confusion can arise,
we call them just trees. 

Given two nodes $v,w\in N(\theta)$, we say that $w \prec v$ if $v$ is
on the path connecting $w$ to the root line. We can identify a line
$\ell$ through the node $v$ it exits by writing $\ell=\ell_{v}$.

We call \textit{internal nodes} the nodes such that there is at least
one line entering them, and \textit{end nodes} the nodes
which have no entering line. We denote by $V(\theta)$
and $E(\theta)$ the set of internal nodes and end nodes,
respectively. Of course $N(\theta)=V(\theta)\cup E(\theta)$.

\begin{figure}[ht]
\vskip-.truecm
\centering 
\ins{92pt}{-85.pt}{$\theta\hskip.2truecm=$}
\includegraphics[width=3in]{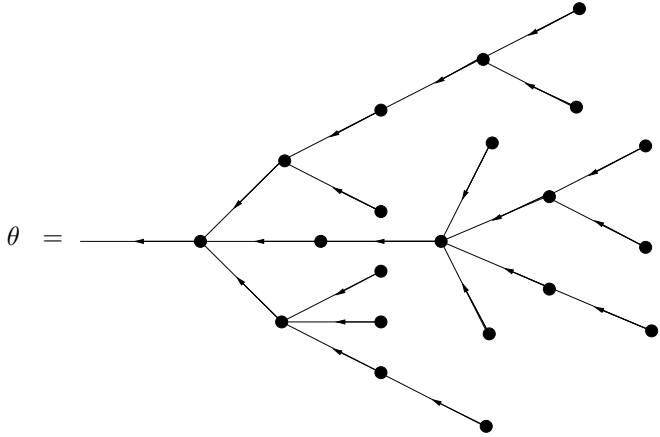}
\caption{An example of unlabelled tree.}
\label{fig:4.1}
\end{figure}

The number of unlabelled trees (i.e. of rooted trees with
no labels) with $N$ nodes -- and hence with $N$ lines --
is bounded by $2^{2N}$, which is a bound on
the number of random walks with $2N$ steps \cite{GM3}.
An example of unlabelled tree is represented in Figure \ref{fig:4.1}.

For each node $v$ denote by $S(v)$ the set of the lines
entering $v$ and set $s_{v}=|S(v)|$; here and henceforth,
given a set $A$, with denote by $|A|$ its cardinality. 
Hence $s_{v}=0$ if $v$ is an end node,
and $s_{v} \geq 1$ if $v$ is an internal node. One has
\be
\sum_{v \in N(\theta)} s_{v} =
\sum_{v \in V(\theta)} s_{v} = |N(\theta)| - 1 ;
\label{eq:4.2} \ee
this can be easily checked by induction
on the number of nodes of the tree.

\subsection{Labels}
\label{sec:4.2}

\noindent
We associate with each node $v \in N(\theta)$ a \textit{mode} label
$\nut_{v}\in\ZZZ^{d}$, and with each line $\ell\in L(\theta)$
a \textit{momentum} label $\nut_{\ell}\in\ZZZ^{d}$, with the
constraints that $\nut_{v} \neq \zerot$ if $v \in E(\theta)$ and

\be
\nut_{\ell_{v}} = \sum_{\substack{ w \in N(\theta) \\ w \preceq v}}
\nut_{w} = \nut_{v}+\sum_{\ell \in S(v)} \nut_{\ell} ,
\label{eq:4.3} \ee
which represents a \textit{conservation rule} for each node.

We also associate with each node $v \in N(\theta)$ an \textit{order} label
$k_{v} \in \{0,1,\ldots,k_{0}\}$, such that
$k_{v}=0$ if $\nut_{\ell_{v}}=\zerot$ and
$k_{v} \ge 1$ if $\nut_{\ell_{v}} \neq \zerot$. We set
\be
k(\theta) = \sum_{v\in N(\theta)} k_{v} , \qquad
\nut(\theta) = \sum_{v \in N(\theta)} \nut_{v} ,
\label{eq:4.4} \ee
which are called the \textit{order} and the \textit{momentum}
of $\theta$, respectively; note that $\nut(\theta)$ is the
momentum of the root line of $\theta$.

Finally we associate with each node $v\in V(\theta)$ a
\textit{badge} label $\rho_{v} \in \{0,1\}$, such that
$k_{v} \in \{0,1\} $ if $\rho_{v}=1$, while $\nut_{\ell_{v}}\neq\zerot$,
$\nut_{v}=0$, and $s_{v}=1$ if $\rho_{v}=0$.

Call $\TT_{k,\nut}$ the set of all trees $\theta$ with order $k$ and
momentum $\nut$, with the constraint that if a line $\ell\in L(\theta)$
has $\nut_{\ell}=0$ and exits a node $v$ with $\nut_{v}=0$ then
$s_{v} \ge 2$. It is easy to check that there exists a
positive constant $\ka$ such that $k(\theta)  \le \ka |N(\theta)|$;
simply use that $s_{v} \neq 1$ when $k_{v}=0$.

\subsection{Diagrammatic rules}
\label{sec:4.3}

\begin{figure}[ht]
\centering 
\ins{173pt}{-005pt}{$=$}
\ins{150pt}{0002pt}{$u^{(k)}_{\nut}$}
\ins{218pt}{-013pt}{$\nut$}
\ins{250pt}{0014pt}{$(k)$}
\includegraphics[width=1in]{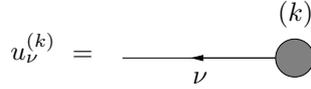}
\caption{Graph element.}
\label{fig:4.2}
\end{figure}

\noindent
We want to show that trees naturally arise when studying
the equations (\ref{eq:2.11}). Let $u^{(k)}_{\nut}$ be represented
with the graph element in Figure \ref{fig:4.2} as a line
with label $\nu$ exiting from a ball with label $(k)$.

\vskip.5truecm

\begin{figure}[ht]
\vskip.4truecm
\centering
\ins{145pt}{-054pt}{$\nut$}
\ins{169pt}{-028pt}{$(k)$}
\ins{205pt}{-048pt}{$=$}
\ins{255pt}{-052pt}{$\nut$}
\ins{270pt}{-052pt}{$\nut_{0}$}
\ins{302pt}{-020pt}{$\nut_{1}$}
\ins{302pt}{0012pt}{$(k_{1})$}
\ins{313pt}{-035pt}{$\nut_{2}$}
\ins{321pt}{-006pt}{$(k_{2})$}
\ins{284pt}{-071pt}{$\nut_{s}$}
\ins{302pt}{-070pt}{$(k_{s})$}
\includegraphics[width=3in]{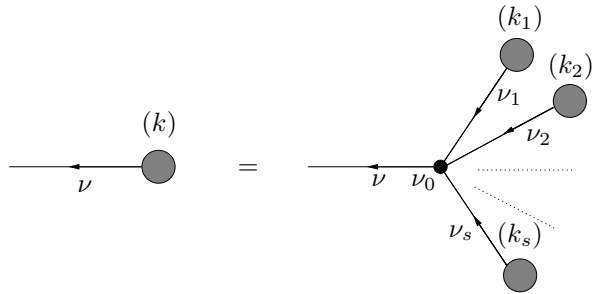}
\caption{Graphical representation of the recursive equations.} 
\label{fig:4.3}
\end{figure}

Then we can represent (\ref{eq:2.11}) graphically as depicted in
Figure \ref{fig:4.3}. Simply expand $[ F(u,\omv t)]^{(k)}_{\nut}$
as in (\ref{eq:2.10}) and represent each factor
$u^{(k_{i})}_{\nut_{i}}$ on the right hand side as a graph element
according to Figure \ref{fig:4.2}. The lines of all such graph elements
enter the same node $v_{0}$. This is a graphical expedient to recall
the conservation rule: the momentum $\nut$ of the root line is the sum
of the mode label $\nut_{0}$ of the node $v_{0}$ plus the sum of the
momenta of the lines entering $v_{0}$. Note that $k_{v_{0}} \ge 1$
as $\nut\neq\zerot$ in (\ref{eq:2.11}).

We represent also (\ref{eq:2.12}) as in Figure \ref{fig:4.3},
with the only difference that now $k_{v_{0}}=0$ and hence
$s_{v_{0}} \ge 2$ (recall the definition of $\TT_{k,\nut}$
at the end of Section \ref{sec:4.2}).

Given any tree $\theta\in\TT_{k,\nut}$ we associate with each
node $v \in N(\theta)$ a \textit{node factor} $\FFF_{v}$ and with each
line $\ell\in L(\theta)$ a \textit{propagator} $\GGG_{\ell}$, by setting
\be
\FFF_{v} := \begin{cases}
(s_{v}!)^{-1} \FF_{s_{v},\nut_{v}} , & \rho_{v} = 1 , \\
- \d_{k_{v}}(\omt\cdot\nu_{\ell_{v}}) \, \id , & \rho_{v} = 0 ,
\end{cases} \qquad\qquad
\GGG_{\ell}:= \begin{cases}
\d_{0}^{-1}(\omt\cdot\nut_{\ell}) \, \id , & \nut_{\ell} \neq \zerot , \\
G , & \nut_{\ell} = \zerot , \end{cases}
\label{eq:4.5} \ee
where $\id$ is the $n \times n$ identity, and $G$ is the $n\times n$
matrix of the form
\be
G = \left(  \begin{matrix} 0 & 0 \\ 0 & -A^{-1} \end{matrix} \right),
\label{eq:4.6} \ee
where the null matrices $0$ are $\tilde n \times \tilde n$,
$\tilde n \times \hat n$, and $\hat n \times \tilde n$, respectively,
while $A$ is the invertible matrix appearing in (\ref{eq:4.1b}).
Define the \textit{value} of the tree $\theta$ as
\be
\Val(\theta) := \Big( \prod_{v\in N(\theta)} \FFF_{v} \Big)
\Big( \prod_{\ell\in L(\theta)} \GGG_{\ell} \Big) .
\label{eq:4.7} \ee
The propagators $\GGG_{\ell}$ are matrices, whereas each $\FFF_{v}$
is a tensor with $s_{v}+1$ indices, which can be associated with
the $s_{v}+1$ lines entering or exiting $v$. In (\ref{eq:4.7})
the indices of the tensors $\FFF_{v}$ must be contracted:
this means that if a node $v$ is connected to a node $v'$ by a line
$\ell$ then the indices of $\FFF_{v}$ and $\FFF_{v'}$ associated with
$\ell$ are equal to each other, and eventually one has to sum
over all the indices except that associated with the root line.

The node factors in (\ref{eq:4.5}) are bounded as
$\max_{j_{1},\ldots,j_{s_{v}+1}}|(\FFF_{v})_{j_{1}\ldots j_{s_{v}+1}}|
\le \Xi_{0}\Xi_{1}^{s_{v}} \rme^{-\xi|\nut_{v}|}$ if $\rho_{v}=1$,
while one has $s_{v}=1$ and $|\d_{0}^{-1}(\omt\cdot\nut_{\ell_{v}})|
\max_{j_{1},j_{2}}|(\FFF_{v})_{j_{1},j_{2}}| \le
|\omt\cdot\nut_{\ell_{v}}|^{(i_{k_{v}}-i_{0})\tau}$ if $\rho_{v}=0$.
As to the propagators one has $\|\GGG_{\ell}\| \le
\g_{0}^{-1} |\nut_{\ell}|^{\tau_{0}}$ for $\nut_{\ell} \neq 0$ and
$\|\GGG_{\ell}\| \le \|A^{-1}\|$ for $\nut_{\ell}=\zerot$,
where $\|\cdot\|$ denotes -- say -- the uniform norm.

By iterating the graphical representation in Figure \ref{fig:4.3}
until only graph elements with $k=1$ appear, one finds
\be
u^{(k)}_{\nut} = \sum_{\theta\in\TT_{k,\nut}} \Val(\theta) ,
\qquad k \geq 1.
\label{eq:4.8} \ee
The \textit{tree expansion} (\ref{eq:4.8}) makes sense since
all node factors and propagators are finite quantities,
and the sum over the labels can be performed. 
Except the mode labels, the last assertion is trivial for all other
labels (as they can assume only a finite number of values).
Finally, the sum over the mode labels is controlled by
the exponential decay $\rme^{-\xi |\nut_{v}|}$
of the Fourier coefficients $\FF_{\nut_{v}}$ of the node factors.

The study of the convergence of the perturbation series is made
difficult by the product of propagators in (\ref{eq:4.7}).
Indeed, the denominators
$\d_{0}(\omt\cdot\nut)$ can be arbitrarily close to zero for
$\nut$ large enough. This problem is usually referred to as
the \textit{small divisor problem}.

\subsection{Compatibility conditions}
\label{sec:4.4}

\noindent
Now, we show that (\ref{eq:4.1}) holds for all the models
considered in Section \ref{sec:3}. Note that to prove (\ref{eq:4.1})
is not a purely technical problem: for the very models of
Section \ref{sec:3}, quasi-periodic solutions in the form of formal
power series or even quasi-periodic solution \textit{tout court}
can fail to exist, if one weaken too much
the assumptions -- see also Section \ref{sec:12.2}.

\begin{figure}[ht]
\vskip-.1truecm
\centering
\ins{010pt}{-027.5pt}{$\theta\;=$}
\ins{194pt}{-027.5pt}{$\theta'\,=$}
\ins{071pt}{-037pt}{$v_{0}$}
\ins{105pt}{-023pt}{$v_{1}$}
\ins{216pt}{-037pt}{$v_{0}$}
\ins{250pt}{-023pt}{$v_{1}$}
\ins{350pt}{-037pt}{$v_{1}$}
\ins{384pt}{-053pt}{$v_{0}$}
\ins{298pt}{-032.5pt}{$=$}
\includegraphics[width=5.4in]{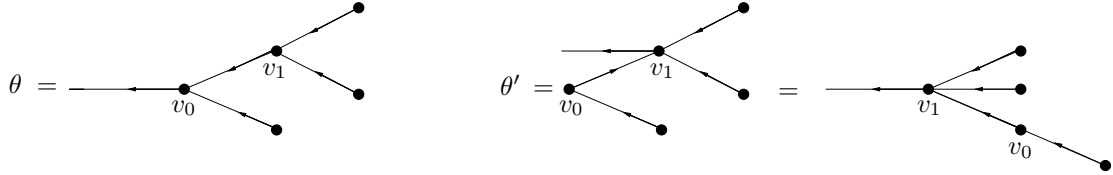}
\caption{An example of tree $\theta'$ obtained from $\theta$ by
detaching the root line from the node $v_{0}$ and
reattaching it to the node $v_{1}$. The arrow of the
line connecting the nodes $v_{0}$ and $v_{1}$ is reverted so
as to point toward the new location of the root. One can always stretch
the lines so as to make all arrows go right to left, as in the
last graph.}
\label{fig:4.4}
\end{figure}

Let us start from the model in Section \ref{sec:3.1};
recall that $d=p=n$ in such a case. We have
already checked that (\ref{eq:3.10}) trivially holds for $k=1$.
Then we can prove by induction that for all $k\ge2$ the
compatibility condition (\ref{eq:3.10}) holds and one can set
$\aaa^{(k)}_{\zerooo}=\zerooo$. The proof proceeds as follows. By
using (\ref{eq:4.8}) we express $\left[ F(u,\omv t) \right]^{(k)}_{\zerot}$
according to (\ref{eq:2.6}) as sum of trees in which all lines
except the root line have nonzero momenta, by the inductive hypothesis.
Given a tree $\theta$ consider together all trees which can be
obtained from $\theta$ by detaching the root line and attaching
it to any other node; see Figure \ref{fig:4.4}.

In that way, we obtain as many trees as nodes of $\theta$;
call $\FF(\theta)$ the set of all such trees.
Of course, all arrows must point toward the root, so that
the trees $\theta'\in\FF(\theta)$ have all the same mode labels
(by construction), but they can have different momenta.
On the other hand, since $\nut=\zerot$, a line
$\ell\in L(\theta')$ either has the same momentum
$\nut_{\ell}$ as in $\theta$ (if the arrow has not been reverted)
or has momentum $-\nut_{\ell}$ (if the arrow has been reverted).
Since $\d_{0}(\nut_{\ell})=\d_{0}(-\nut_{\ell})$, this means that
the corresponding propagator $\GGG_{\ell}$ does not change.
The combinatorial factors of the trees $\theta'\in\FF(\theta)$
are in general different from those of $\theta$ (because
the values of $s_{v}$, $v\in N(\theta')$, can change),
but if we sum together all nonequivalent trees we realise
(with a little effort: one must perform the computation
to convince himself that the assertion is true!)
that we obtain a common value times a factor $\rmi\nnn_{v_{0}}$,
if $v_{0}$ is the node which the root line is attached to.
Therefore, since
\be
\nut = \sum_{v_{0}\in N(\theta)} \nut_{v_{0}} = 0 ,
\label{eq:4.9} \ee
the sum of all the tree values gives zero.
This implies (\ref{eq:4.1}) with $\tilde n = n$ and $\hat n =0$,
and the coefficient $\aaa^{(k)}_{\zerooo}$ is left arbitrary,
and it can be chosen to be zero.

In fact, the argument above does not depend on the value of the
coefficients $\aaa^{(k)}_{\zerooo}$, which therefore can be arbitrarily
chosen; in particular they can and will arbitrarily fixed to be zero.
This was expected: changing $\aaa^{(k)}_{\zerooo}$ means changing
the constant $\aaa_{0}$ in (\ref{eq:3.4}), which is arbitrary
since it is the vector parameterising the torus.

The case of the standard map -- see Section \ref{sec:3.2} --
can be discussed in the same way. We omit the details.

In the case of the model in Section \ref{sec:3.3}, the identity
$\left[ \partial_{\aaa} f(\aaa,\bbb) \right]^{(k)}_{\zerooo}=\zerooo$
for $k \ge 2$ can be proved as above, by relying on the same
cancellation mechanism. The compatibility condition
$\left[ \partial_{\bbb} f(\aaa,\bbb) \right]^{(k)}_{\zerooo}=\zerooo$
for $k\ge 2$ can be imposed by using (\ref{eq:3.25}) and
fixing $\bbb^{(k)}_{\zerooo}=-
( \partial_{\bbb}^{2} f_{\zerooo}(\bbb_{0}) )^{-1}
\boldsymbol{\Phi}_{k}$. This implies once more (\ref{eq:4.1})
with $\tilde n=r$ and $\hat n=s$. Again the coefficients
$\aaa^{(k)}_{\zerooo}$ can be arbitrarily set to be zero.

Finally for the model in Section \ref{sec:3.4} one can
use (\ref{eq:3.31}) to obtain (\ref{eq:4.1}) with $\tilde n=0$
and $\hat n = n =1$.

\zerarcounters
\section{Multiscale analysis}
\label{sec:5}

\noindent
To be able to bound the tree value (\ref{eq:4.7}), we need
to control the product of propagators. This will be done through 
a \textit{multiscale analysis}. To this aim, for each tree line
we introduce a new label characterising the size of the corresponding
propagator, that we call the \textit{scale} label.

Essentially, we say that $\nu\in\ZZZ^{d} \setminus \{\zerot\}$ is on scale
\be
\begin{cases}
n \geq 1, & \hbox{if } 2^{-n} \g \leq |\omt\cdot\nut| <
2^{-(n-1)}\g , \\
n = 0 , & \hbox{if } \g \leq |\omt\cdot\nut| , \end{cases}
\label{eq:5.1} \ee
where $\g$ is the constant appearing in (\ref{eq:1.5}),
and we say that a line $\ell$ has a scale label $n_{\ell}=n$ if
$\nu_{\ell}$ is on scale $n$.

As a matter of fact, in practice the sharp multiscale decomposition
in (\ref{eq:5.1}) is a little annoying because, as we shall see,
we have to consider derivatives. Thus, it is actually more convenient
to replace it with a smooth decomposition through $C^{\io}$
compact support functions. Let $\psi$ be a nondecreasing $C^{\infty}$
function defined in $\RRR_{+}$, such that
\be
\psi(x) = \left\{
\begin{array}{ll}
1 , & \text{for } x \geq \g , \\
0 , & \text{for } x \leq \g/2 ,
\end{array} \right.
\label{eq:5.2} \ee
and set $\chi(x) := 1-\psi(x)$. For all $n \in \ZZZ_{+} = \NNN
\cup \{0\}$ define $\chi_{n}(x) := \chi(2^{n}x)$
and $\psi_{n}(x) := \psi(2^{n}x)$, and set
\be
\Xi_{n}(x)=\chi_{0}(|x|)\ldots \chi_{n-1}(|x|) \chi_{n}(|x|) , \qquad
\Psi_{n}(x)=\chi_{0}(|x|)\ldots \chi_{n-1}(|x|) \psi_{n}(|x|) ,
\label{eq:5.3} \ee
where $\Psi_{0}(x)$ is meant as $\Psi_{0}(x)=\psi_{0}(|x|)$.

Then we change the definition of the propagator to be
associated with each $\ell$  with $\nut_{\ell}\neq\zerot$, by
associating with each such line $\ell$ a scale label
$n_{\ell} \in \ZZZ_{+}$ and a propagator
\be
\GGG_{\ell} = G^{[n_{\ell}]}(\omt\cdot\nut_{\ell}) , \qquad
G^{[n]}(\omt\cdot\nut) := \Psi_{n}(\omt\cdot\nut)\,
\d_{0}^{-1}(\omt\cdot\nut) \, \id ,
\label{eq:5.4} \ee
which replaces the previous definition in (\ref{eq:4.5}). If
$\Psi_{n}(x) \neq 0$ then $2^{-n-1}\g \le |x| \le 2^{-n+1}\g$, so that
for any each $\ell \in L(\theta)$ one has either $\GGG_{\ell}=0$
or $\|\GGG_{\ell}\| \le \g_{0}^{-1}2^{(n_{\ell}+1)i_{0}}$.
For completeness we also associate a scale label $n_{\ell}=-1$
with each line $\ell$ with momentum $\nut_{\ell}=\zerot$.
Note that, while with the sharp decomposition (\ref{eq:5.1}) a
momentum $\nut$ identifies uniquely the scale $n$, on the contrary
by using the smooth decomposition for each momentum $\nut$
there are two possible (adjacent) values $n$ such that
$G^{[n]}(\omt\cdot\nut)\neq0$.

A tree expansion like (\ref{eq:4.8}) still holds, with the
difference that now the trees $\theta\in\TT_{k,\nut}$ carry also
the scale labels, and we have to sum also on these labels.
The equality between the two expansions follows immediately from
the observation that $\sum_{n=0}^{\io} \Psi_{n}(x)=1$ for
all $x\in\RRR \setminus\{0\}$.

If $\gotN_{n}(\theta)$ denotes the number of lines $\ell\in L(\theta)$
with scale $n_{\ell}=n$, then we can bound in (\ref{eq:4.7})
\be
\prod_{\ell\in L(\theta)} \left\| \GGG_{\ell} \right\| \leq
\g_{0}^{-k} 2^{k i_{0}} \prod_{n=0}^{\io} 2^{n i_{0} \gotN_{n}(\theta)} , 
\label{eq:5.5} \ee
with $i_{0}$ defined in (\ref{eq:1.2}),
so that the problem is reduced to bounding $\gotN_{n}(\theta)$.

The product of propagators gives problems when the small
divisors ``accumulate''. To make more precise the idea of
accumulation we introduce the notion of cluster.
Once all lines of a tree $\theta$ have been given their scale labels,
for any $n\geq 0$ we can identify the maximal connected sets of lines
with scale not larger than $n$. If at least one among such lines
has scale equal to $n$ we say that the set is a \textit{cluster}
on scale $n$. Given a cluster $T$
call $L(T)$ the set of lines of $\theta$ contained in $T$,
and denote by $N(T)$ the set of nodes connected by such lines.
We define $k(T)=\sum_{v\in N(T)}k_{v}$
the \textit{order of the cluster} $T$.

Any cluster has either one or no exiting line,
and can have an arbitrary number of entering lines.
We call \textit{self-energy clusters} the clusters which have
one exiting line and only one entering line and are such that both
lines have the same momentum; the terminology is borrowed from
quantum field theory. This means that if $T$
is a self-energy cluster and $\ell_{1}$ and $\ell_{2}$ are the
lines entering and exiting $T$, respectively, then
$\nu_{\ell_{1}}=\nu_{\ell_{2}}$, so that
\be
\sum_{v\in N(T)} \nut_{v} = 0 .
\label{eq:5.6} \ee
By construction the scales of the lines $\ell_{1}$
and $\ell_{2}$ can differ at most by 1, and setting
$n_{T}=\min\{n_{\ell_{1}},n_{\ell_{2}}\}$, by definition of cluster
one has $n_{\ell} < n_{T}$ for all $\ell\in L(T)$.

We define the \textit{value} of the self-energy cluster $T$ whose
entering line has momentum $\nu$ as the matrix
\be
\VV_{T}(\omt\cdot\nut) := \Big( \prod_{v\in N(T)} \FFF_{v} \Big)
\Big( \prod_{\ell\in L(T)} \GGG_{\ell} \Big) ,
\label{eq:5.7} \ee
where all the indices of the node factors must be contracted
except those associated with the line $\ell_{1}$
entering $T$ and with the line $\ell_{2}$ exiting $T$.

We can extend the notion of self-energy cluster also to
a single node, by saying that $v$ is a self-energy cluster
if $s_{v}=1$ and the line entering $v$ has the same momentum
as the exiting line. In that case (\ref{eq:5.7}) has to be
interpreted as $\VV_{T}(\omt\cdot\nut)=\FFF_{v}$:
in particular it is independent of $\omt\cdot\nut$ if $\rho_{v}=1$. 
If $T$ consists of only one node (and hence contains no line)
we say that $T$ is a cluster on scale $-1$.

The simplest self-energy cluster one can think of consists of
only one node $v$, but then (\ref{eq:5.6}) implies $\nut_{v}=0$.
For the model in Section \ref{sec:3.1}, the corresponding
value is zero, and hence the simplest nontrivial self-energy clusters
contain at least two nodes. On the contrary for the models
in Sections \ref{sec:3.3} and \ref{sec:3.4} one can also
have clusters with only one node.

The reason why it is important to introduce the self-energy clusters
is that if we could neglect them then the product of small divisors
would be controlled. Indeed, let us denote by $\gotR_{n}(\theta)$
the number of lines on scale $n$ which do exit a self-energy cluster,
and set $\gotN_{n}^{*}(\theta)=\gotN_{n}(\theta)-\gotR_{n}(\theta)$. Then
an important result, known as the \textit{Siegel-Bryuno lemma}, is that
\be
\gotN_{n}^{*}(\theta) \leq c\,2^{-n/\tau} K(\theta) ,
\qquad K(\theta) := \sum_{v \in N(\theta)} |\nut_{v}| ,
\label{eq:5.8} \ee
for some constant $c$, where $\tau$ is the Diophantine exponent
in (\ref{eq:1.5}). See Appendix \ref{app:A} for a proof.

If no self-energy clusters could occur (so that $\gotR_{n}(\theta)=0$)
the Siegel-Bryuno lemma would allow us to bound in (\ref{eq:5.5})
\be
\prod_{n=0}^{\io} 2^{n i_{0} \gotN_{n}(\theta)} =
\prod_{n=0}^{\io} 2^{n i_{0} \gotN_{n}^{*}(\theta)} \leq
2^{n_{0}i_{0}k} \prod_{n=n_{0}+1}^{\io}
2^{n i_{0} \gotN_{n}^{*}(\theta)} \leq
C_{1}^{k} \exp \Big( \xi(n_{0}) K(\theta) \Big) , 
\label{eq:5.9} \ee
with $C_{1}=2^{n_{0}i_{0}}$ and
\be \xi(n_{0}) := i_{0} c \sum_{n=n_{0}+1}^{\io} n 2^{-n/\tau} .
\label{eq:5.10} \ee
Since $\xi(n_{0}) \to 0$ as $n_{0}\to\io$ one can fix $n_{0}$
in such a way that $\xi(n_{0}) \le \xi/4$ (the constant $\xi$
being defined after (\ref{eq:2.6})). Then, by extracting a factor
$\rme^{-\xi|\nut_{v}|/2}$ from each node factor $\FFF_{v}$,
one could easily perform the sum over the Fourier labels,
so as to obtain an overall bound $C_{2}^{k}\rme^{-\xi|\nut|/2}$
on $u^{(k)}_{\nut}$ for a suitable constant $C_{2}$. This would imply
the convergence of the perturbation series (\ref{eq:2.1})
for $\e$ small enough, say for $|\e|<\e_{0}$ for some $\e_{0}>0$.
However, there are self-energy clusters and they produce factorials,
as the example in Appendix \ref{app:A} shows,
so that we have to deal with them.

\zerarcounters
\section{Resummation of the series}
\label{sec:6}

\noindent
Let us come back to the equation (\ref{eq:2.11}).
To simplify the analysis, let us initially assume that we are using
the sharp multiscale decomposition (\ref{eq:5.1}), so that
each momentum fixes uniquely the corresponding scale.
If we take the tree expansion of the right hand side
of (\ref{eq:2.11}), according to the
diagrammatic rules described in Section \ref{sec:4},
we can distinguish between contributions in which the root line exits
a self-energy cluster $T$, that we can write as
\be
\sum_{T:k(T)<k} \VV_{T}(\omt\cdot\nut) \, u^{(k-k(T))}_{\nut} ,
\label{eq:6.1} \ee
and all the other contributions, that we denote by
$[F(u,\omv t)]^{(k-1)*}_{\nut}$. In (\ref{eq:6.1}) both
the entering and exiting lines of $T$ have the same scale
$n_{T}$, and the sum is over all clusters $T$ on scale $<n_{T}$.

By writing (\ref{eq:2.11}) as
\be
\d_{0}(\omt\cdot\nut) \, u^{(k)}_{\nut} = - 
\sum_{p=1}^{\min\{k,k_{0}\}} \d_{p}(\omt\cdot\nut) u^{(k-p)}_{\nut} +
\left[ F(u,\omv t) \right]^{(k-1)}_{\nut} 
\label{eq:6.2} \ee
we can shift the contributions (\ref{eq:6.1}) to the left hand side
of (\ref{eq:6.2}), so as to obtain
\be
\d_{0}(\omt\cdot\nut) \, u^{(k)}_{\nut} - \sum_{T : k(T)<k}
\VV_{T}(\omt\cdot\nut)\,u^{(k-k(T))}_{\nut} = 
\left[ F(u,\omv t) \right]^{(k-1)*}_{\nut} .
\label{eq:6.3} \ee
By summing over $k$ and setting
\be
M(\omt\cdot\nut,\e) : = \sum_{k=1}^{\io} \e^{k}
\sum_{T:k(T)=k}\VV_{T}(\omt\cdot\nut) ,
\label{eq:6.4} \ee
then (\ref{eq:6.3}) gives, formally,
\be
\DD(\omt\cdot\nut,\e)\,u_{\nu} =
\left[ F(u,\omv t) \right]^{*}_{\nut} .
\qquad \DD(\omt\cdot\nut,\e) :=
\d_{0}(\omt\cdot\nut) \, \id - M(\omt\cdot\nut,\e) .
\label{eq:6.5} \ee
The motivation for proceeding in this way is that, at the price
of changing $\d_{0}(\omt\cdot\nut)$ into $\DD(\omt\cdot\nut,\e)$, hence of
changing the propagators, lines exiting self-energy clusters
no longer appear. Therefore, in the tree expansion of
the right hand side of the equation, we have eliminated
the self-energy clusters, that is the source of the problem
of accumulation of small divisors.

Unfortunately the procedure described above has a problem:
$M(\omt\cdot\nut,\e)$ itself is a sum of self-energy clusters, which
can still contain some other self-energy clusters on lower scales.
So finding a good bound for $M(\omt\cdot\nut,\e)$ could have
the same problems as for the values of the trees.

To deal with such a difficulty we modify the prescription
by proceeding recursively, in the following sense.
Let us start from the momenta $\nut$ which are on scale $n=0$.
Since the only self-energy clusters $T$ with $n_{T}=0$
are those (on scale $-1$) containing only one node,
for such $\nut$ the matrix $M(\omt\cdot\nut,\e)$ is just a constant.
Next, we pass to the momenta $\nut$ which are on scale $n=1$,
and we consider (\ref{eq:6.5}) for such $\nut$;
now all self-energy clusters $T$ whose values contribute to
$M(\omt\cdot\nut,\e)$ cannot contain any self-energy clusters,
because the lines $\ell\in L(T)$ are on scale $n_{\ell}=0$.
Then, we consider the momenta $\nut$ which are on scale $n=2$:
again all the self-energy clusters contributing to
$M(\omt\cdot\nut,\e)$ do not contain any self-energy clusters,
because the lines on scale $n=0,1$ cannot exit
self-energy clusters by the construction of the previous step,
and so on. The conclusion is that we have obtained a different
expansion for $u(\omt t,\e)$, that we call a \textit{resummed series},
\be
u(\omt t,\e) = \sum_{\nut\in\ZZZ^{d}} \rme^{\rmi\omt\cdot\nut t}
u_{\nut} , \qquad u_{\nut} =
\sum_{k=1}^{\io} \e^{k} u^{[k]}_{\nut}(\e) ,
\label{eq:6.6} \ee
where the self-energy clusters do not appear any more
in the tree expansion and the propagators must be defined recursively,
as follows. The propagator $\GGG_{\ell}$ of a line $\ell$
on scale $n_{\ell}=n$ and momentum $\nut_{\ell}=\nut$ is the matrix
\be
\GGG_{\ell} := G^{[n]}(\omt\cdot\nut,\e) = \left( \d_{0}(\omt\cdot\nut)
\, \id - \MM^{[n-1]}(\omt\cdot\nut,\e) \right)^{-1} ,
\label{eq:6.7} \ee
with
\be
\MM^{[n]}(\omt\cdot\nut,\e) := \sum_{T {\rm on \, scale \,} \le n}
\e^{k(T)} \, \VV_{T}(\omt\cdot\nut) , 
\label{eq:6.8} \ee
where the value $\VV_{T}(\omt\cdot\nut)$ is written in accord
with (\ref{eq:5.7}), with all the lines $\ell'\in L(T)$
on scales $n_{\ell'}<n$ and the corresponding propagators
$\GGG_{\ell'}$ expressed in terms of matrices
$\MM^{[n_{\ell'}]}(\omt\cdot\nut_{\ell'},\e)$ as in (\ref{eq:6.7}).

By construction, the new propagators depend on $\e$,
so that the coefficients $u^{[k]}_{\nut}(\e)$ depend explicitly
on $\e$: hence (\ref{eq:6.6}) is not a power series expansion.

If we use the smooth multiscale decomposition, then the
algorithm above must be suitably modified. We define 
recursively the propagators $\GGG_{\ell}=G^{[n_{\ell}]}
(\omt\cdot\nut_{\ell},\e)$ by setting for $n\ge0$
\begin{subequations}
\begin{align}
& G^{[n]}(\omt\cdot\nut,\e) =
\Psi_{n}(\omt\cdot\nut) \left( \d_{0}(\omt\cdot\nut) \id -
\MM^{[n-1]}(\omt\cdot\nut,\e) \right)^{-1}
\label{eq:6.9a} \\
& \MM^{[n]}(x,\e) = \MM^{[n-1]}(x,\e) +
\Xi_{n}(x) \, M^{[n]}(x,\e) , \qquad
M^{[n]}(x,\e) = \!\! \sum_{T\in\RR_{n}} \e^{k(T)} \, \VV_{T}(x) ,
\label{eq:6.9b}
\end{align}
\label{eq:6.9} \end{subequations}
\vskip-.3truecm
\noindent where $\RR_{n}$ is the set of self-energy clusters
on scale $n$ which do not contain any other self-energy clusters.

The matrices $\MM^{[n]}(\omt\cdot\nut,\e)$ are called the
\textit{self-energies}.
The new propagators (\ref{eq:6.9a}) are called, by exploiting once
more the analogy with quantum field theory,
the \textit{dressed propagators}.

The coefficients $u^{[k]}_{\nu}(\e)$ still admit a tree expansion
\be
u^{[k]}_{\nu}(\e) =
\sum_{\theta\in\TT^{\RR}_{k,\nut}} \Val(\theta) , \qquad
\Val(\theta) := \Big(\prod_{v\in N(\theta)}\FFF_{v}\Big)
\Big(\prod_{\ell\in L(\theta)} \GGG_{\ell} \Big)
\qquad \nut \neq \zerot, \quad k\geq1 ,
\label{eq:6.10} \ee
which replaces (\ref{eq:4.8}). In particular
$\TT^{\RR}_{k,\nu}$ is defined as the set of \textit{renormalised
trees} of order $k$ and momentum $\nu$, where ``renormalised''
means that the trees do not contain any self-energy clusters.

Since for any tree $\theta\in\TT^{\RR}_{k,\nu}$
one has $\gotN_{n}(\theta)=\gotN_{n}^{*}(\theta)$,
we can bound the product of propagators according to (\ref{eq:5.5})
and (\ref{eq:5.9}), provided the propagators on scale $n$
can still be bounded proportionally to $2^{n i_{0}}$.
In general there is no reason why this should occur,
because of the extra term $\MM^{[n-1]}(\omt\cdot\nut,\e)$
appearing in (\ref{eq:6.9a}).

The discussion of such an issue depends on the particular model
one is studying. We shall see in the next sections what happens
for the models considered in Section \ref{sec:3}.
We shall first consider cases in which the dressed propagators
can be essentially bounded as the old ones, and then cases
in which this is no longer true. By modifying further
the resummation procedure described above, we shall see that
something can still be achieved also in these cases.

\zerarcounters
\section{Cancellations and convergence of the series -- maximal tori}
\label{sec:7}

\noindent
Let us consider the matrix $\MM^{[n]}(x,\e)$ introduced
in (\ref{eq:6.9b}), and let us study its dependence on the
first argument $x=\omt\cdot\nut$ for the models of Section \ref{sec:3}.
As usual, let us consider first the model in Section \ref{sec:3.1}.

It is a remarkable cancellation that $\MM^{[n]}(x,\e)$ vanishes
in $x$ up to second order, that is $\MM^{[n]}(x,\e)=O(x^{2})$.
The symmetry properties
\be
\MM^{[n]}(x,\e)=(\MM^{[n]}(-x,\e))^{T}=(\MM^{[n]}(x,\e))^{\dagger} ,
\label{eq:7.1} \ee
with $T$ and $\dagger$ denoting transposition and adjointness,
are essential for the proof (such properties are trivially satisfied
for $n=0$, and can be proved by induction on $n$ -- see \cite{GM2,GBG}
for more details). Indeed, by using (\ref{eq:7.1}),
the cancellation $\MM^{[n]}(0,\e) = 0$ can be proved as
the cancellation $\left[ F(u,\omv t) \right]^{(k)}_{\zerot}=
\left[ \partial_{\aaa} f(\aaa) \right]^{(k)}_{\zerooo}=0$
discussed in Section \ref{sec:4.4}, with the exiting line of
the self-energy clusters playing the role of the root line.
The first order cancellation requires $\partial_{x}\MM^{[n]}(0,\e)=0$,
and this can be proved through a similar cancellation mechanism:
besides the exiting line one has to detach also the entering line
and reattach it to all the other nodes inside the
self-energy clusters. If the function $f$ in (\ref{eq:3.1})
is even in $\aaa$, then the first order cancellation
follows also from parity properties \cite{Ga1}.

Both the cancellations and the symmetry properties are only formal
as far as we have not proved that the self-energies are
well-defined quantities. To this aim we need to control the product
of propagators in (\ref{eq:5.7}), with the propagators defined
according to (\ref{eq:6.9}). An important ingredient of the analysis
is that also for all self-energy clusters one can prove a
bound like (\ref{eq:5.8}). More precisely, if we denote
by $\gotN_{n}(T)$ the number of lines $\ell\in L(T)$
on scale $n$, for all $n\ge0$ and all $T\in \RR_{n}$ one has
\be
\gotN_{n'}^{*}(T) \leq c'\,2^{-n'/\tau} K(T) ,
\qquad K(T) := \sum_{v \in N(T)} |\nut_{v}| , \qquad n'\le n , 
\label{eq:7.2} \ee
for some constant $c'$. To prove (\ref{eq:7.2}) one first
show that, for all $n\ge0$ and all $T\in \RR_{n}$, one has
\be
\sum_{v\in N(T)} |\nut_{v}| > c'' 2^{n/\tau} ,
\label{eq:7.3} \ee
for some constant $c''$, then one proceeds by induction on
the order $k(T)$; see Appendix \ref{app:B} for details.

Therefore, if we were able to prove for the dressed propagators
an estimate like $\|G^{[n]}(x,\e)\|\le 2/x^{2}$,
then we could use (\ref{eq:7.2}) and (\ref{eq:7.3}) to bound
\be
\prod_{\ell \in L(T)} \left\| \GGG_{\ell} \right\| \le
\g_{0}^{-k} 2^{3k} \prod_{n'=0}^{n} 2^{2n'\gotN_{n'}^{*}(T)} ,
\label{eq:7.4} \ee
in such a way to obtain
\be
\left| \VV_{T}(\omt\cdot\nut) \right| \le C_{2}^{k(T)} \rme^{-\xi K(T)/2} ,
\label{eq:7.5} \ee
for a suitable constant $C_{2}$, independent of $T$. In particular,
this would ensure the well-definedness of the self-energies.
At this point, the cancellations would allow us to 
write, for some constant $C$ and for all $n\ge 0$,
\be
\MM^{[n]}(x,\e) = \e^{2}x^{2}\overline\MM^{[n]}(x,\e) ,
\qquad \left\| \overline\MM^{[n]}(x,\e) \right\| \leq C ,
\label{eq:7.6} \ee
where we have taken into account also that $\VV_{T}(x) \neq 0$
requires $k(T) \geq 2$ (cf. Section \ref{sec:5}).
In turn (\ref{eq:7.6}) would implies that, for $\e$ small enough,
\be
\left\| G^{[n]}(x,\e) \right\| = \left\|
\left( x^{2} - \e^{2}x^{2}
\overline\MM^{[n-1]}(x,\e) \right)^{-1} \right\|
\leq \frac{2}{x^{2}} .
\label{eq:7.7} \ee

The bounds (\ref{eq:7.7}) on the propagators, the symmetry properties,
and the cancellations are proved all together, as follows.
First note that, if the second order cancellation holds, one can write
\be
\e^{2} \overline\MM^{[n]}(x,\e) = \int_{0}^{1} {\rm d}t
\, \left( 1 - t \right) \partial_{x}^{2} \MM^{[n]}(t x,\e) ,
\label{eq:7.8} \ee
so that the bound in (\ref{eq:7.6}) is essentially a bound on the
second derivative of the self-energies. The case $n=0$ is easily
checked. Then the proof proceeds by induction, by relying on the
recursive definition of $\MM^{[n]}(x,\e)$ -- see (\ref{eq:6.9})
and (\ref{eq:5.7}) --, and taking advantage of the smooth
multiscale decomposition to perform the derivatives.
More precisely, we assume that both the cancellations -- and hence the
bounds (\ref{eq:7.6}) -- and the symmetry properties (\ref{eq:7.1})
hold for all $n'<n$. This means that all the dressed propagators
of the lines on scales $\le n$ are bounded according to (\ref{eq:7.7}),
so that we can use the bound (\ref{eq:7.5}) to prove that
also $\MM^{[n]}(x,\e)$ is well-defined. Then the cancellation
mechanism described at the beginning of the section shows
that also at the step $n$ the symmetry properties (\ref{eq:7.1})
and the cancellations leading to (\ref{eq:7.6}) are satisfied;
in particular also the propagators $G^{[n+1]}(x,\e)$ are
bounded proportionally to $|x|^{2}$ according to (\ref{eq:7.7}).

The conclusion is that the series in (\ref{eq:6.6}) for
$u_{\nut}=\aaa_{\nnn}$ converges for $\e$ small enough.
Therefore, the function $u(\omt t,\e)=\aaa(\ooo t,\e)$ is analytic
in $\e$ (notwithstanding that the expansion in $\e$ is not a
power expansion), so that we can say \textit{a posteriori}
that the original power series (\ref{eq:3.14}) also converges.
It is straightforward to see that $\aaa_{\nnn}$ decays
exponentially in $\nnn$, which implies that the function
$\aaa(\boldsymbol{\psi},\e)$ is also analytic in $\boldsymbol{\psi}$.

The case of the standard map -- see Section \ref{sec:3.2} -- can be
discussed in the same way. We do not repeat the analysis and
refer to \cite{G1,BeG2} for details.

\zerarcounters
\section{Summation of the divergent series -- hyperbolic tori}
\label{sec:8}

\noindent
Now we consider the model introduced in Section \ref{sec:3.3}.
In that case, one has
\be
\MM^{[n]}(x,\e) = \left( \begin{matrix}
0 & 0 \\ 0 & \e \partial_{\bbb}^{2}f_{\zerooo}(\bbb_{0})
\end{matrix} \right) + O(\e^{2}) ,
\label{eq:8.1} \ee
so that, already only keeping the first order terms, one realises
that a cancellation like (\ref{eq:7.6}) cannot expected to hold.
Indeed, in order to study the convergence of the series (\ref{eq:6.6}),
we need at least the perturbation series (\ref{eq:2.1}) to be
formally well-defined to all orders; in turn this requires
matrix $\partial_{\bbb}^{2}f_{\zerooo}(\bbb_{0})$ to be nonsingular
-- see Sections \ref{sec:3.4} and \ref{sec:4.4} --
and hence different from $0$.

Let us assume first that the matrix
$\partial_{\bbb}^{2}f_{\zerooo}(\bbb_{0})$ is positive definite,
that is that its eigenvalues $a_{1},\ldots,a_{s}$ are positive,
i.e $a_{i}>0$ for $i=1,\ldots,s$ (in particular this means that
$\bbb_{0}$ is a maximum point for the function $f_{\zerooo}(\bbb)$).

We write both the self-energies and the propagators as block matrices,
\be
\MM^{[n]}(x,\e) = \left( \begin{matrix}
\MM^{[n]}_{\aaa\aaa}(x,\e) & \MM^{[n]}_{\aaa\bbb}(x,\e) \\
\MM^{[n]}_{\bbb\aaa}(x,\e) & \MM^{[n]}_{\bbb\bbb}(x,\e)
\end{matrix} \right), \qquad
G^{[n]}(x,\e) = \left( \begin{matrix}
G^{[n]}_{\aaa\aaa}(x,\e) & G^{[n]}_{\aaa\bbb}(x,\e) \\
G^{[n]}_{\bbb\aaa}(x,\e) & G^{[n]}_{\bbb\bbb}(x,\e)
\end{matrix} \right),
\label{eq:8.2} \ee
where the four blocks are $r \times r$, $r \times s$, 
$s \times r$, and $s \times s$ matrices, respectively.

Then, one can prove that the parity properties (\ref{eq:7.1})
still hold, and moreover, formally, one has the cancellations
\be
\MM^{[n]}_{\aaa\aaa}(x,\e) = O(\e^{2}x^{2}) , \qquad
\MM^{[n]}_{\aaa\bbb}(x,\e) = \MM^{[n]}_{\bbb\aaa}(x,\e) = O(\e^{2}x) .
\label{eq:8.3} \ee
The proof of such assertions can be performed by induction,
and follows the same pattern as described in
Section \ref{sec:7} -- see \cite{GG1,GBG} for details.
We have used the word ``formally'' because the cancellations hold
as far as the dressed propagators (\ref{eq:6.9a})
can be bounded essentially as the old ones (\ref{eq:5.4}) -- a
property that we have not yet proved.

The main implication of (\ref{eq:7.1}) and (\ref{eq:8.3}) is that
the eigenvalues $\l^{[n]}_{i}(x,\e)$ of the self-energies
$\MM^{[n]}(x,\e)$ are of the form
\be
\l^{[n]}_{i}(x,\e) = \begin{cases}
O(\e^{2}x^{2}) , & \quad i =1,\ldots,r , \\
a_{i-r}\e + O(\e^{2}) , & \quad i =r+1,\ldots,d . \end{cases}
\label{eq:8.4} \ee
In particular, if $\e<0$, the eigenvalues $x^{2} - \l^{[n]}_{i}(x,\e)$
of the matrices $\d_{0}(\omt\cdot\nut)\,\id - \MM^{[n]}(x,\e)$
are such that $x^{2} - \l^{[n]}_{i}(x,\e) \ge x^{2}/2$ for
$i=1,\ldots,r$ and $x^{2} - \l^{[n]}_{i}(x,\e) \ge x^{2} +
| \e a_{i-r}|/2$ for $i=r+1,\ldots,d$,
provided $\e$ is small enough and the higher order
corrections in (\ref{eq:8.4}) remain small. The last property
is automatically satisfied if the block matrices in (\ref{eq:8.3})
are dominated by the first nontrivial orders.

All the properties described above become rigorous if
the dressed propagators $G^{[n]}(x,\e)$ are bounded proportionally
to $x^{-2}$ (more generally any power of $|x|$ would suit), say
\be
\left\| G^{[n]}(x,\e) \right\| \le \frac{2}{x^{2}} .
\label{eq:8.5} \ee
So, all we have to do is to prove together all the above properties
(\ref{eq:7.1}) and (\ref{eq:8.3}), by induction.
Indeed, for $n=0$ the properties are trivially satisfied, and
at any step $n$, by the inductive hypothesis, the bounds
$\| G^{[n']}(x,\e) \| \le 2 / x^{2}$ are satisfied for all
$n' \le n$, so that both (\ref{eq:7.1}) and (\ref{eq:8.3})
hold for $n$, and in turn this implies (\ref{eq:8.4}) and hence
the bound (\ref{eq:8.5}) for $n+1$.

\begin{figure}[ht]
\vskip-.1truecm
\centering
\ins{284pt}{-60pt}{$\hbox{Re}\,\e$}
\ins{224pt}{-00pt}{$\hbox{Im}\,\e$}
\includegraphics[width=2in]{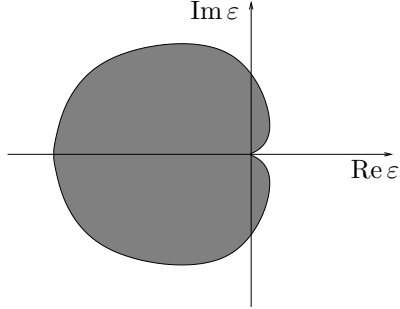}
\caption{Analyticity domain for the hyperbolic invariant torus.}
\label{fig:8.1}
\end{figure}

Therefore, the series (\ref{eq:6.6}) converges for $\e$ small enough,
even if analyticity is prevented because of the condition $\e<0$;
then we say that the perturbation series (\ref{eq:2.1})
is a \textit{divergent series} (in point of fact, it is very likely
that it does not converge, though there is no proof of that).
The function is only $C^{\io}$ in $\e$ at $\e=0$
-- such a result improves a previous one by Treshch\"{e}v \cite{T},
where $C^{\io}$-smoothness in $\sqrt{\e}$ was proved at $\e=0$.
In fact, one can say a little more about the dependence of the
invariant torus on the perturbation parameter $\e$ in the
complex $\e$-plane: the lower-dimensional torus turns out to be
analytic in the heart-shaped domain of Figure \ref{fig:8.1} \cite{GG1}.

As said in Section \ref{sec:3.3} the torus will be said hyperbolic
in that case. This is a somewhat improper terminology. Indeed,
when studying lower-dimensional tori one usually considers
Hamiltonians of the form \cite{Gr,Me2,E1,Poe2}
\be
\calH(\aaa,\qqq,\AAA,\ppp) = \frac{1}{2} \AAA^{2} +
\frac{1}{2} \ppp^{2} + \sum_{i=1}^{s} \l_{i} q_{i}^{2} +
\e f(\aaa,\qqq) ,
\label{eq:8.6} \ee
(or generalisations of its), with $(\aaa,\AAA)\in\TTT^{r}\times\RRR^{r}$
and $(\qqq,\ppp)\in\RRR^{s}\times\RRR^{s}$. Thus, for $\e=0$
the coordinates $\aaa$ freely rotates with some
rotation vector $\ooo=(\o_{1},\ldots,\o_{r})$, while the
coordinates $(\qqq,\ppp)$ moves around an equilibrium point,
which is elliptic if $\l_{i}>0$ for all $i=1,\ldots,s$
and hyperbolic if $\l_{i}<0$ for all $i=1,\ldots,s$.
The numbers $\o_{1},\ldots,\o_{r}$ are called the
\textit{proper frequencies}, while the numbers $\l_{1},\ldots,\l_{s}$
are called the \textit{normal frequencies}.

Then one can study the problem of persistence of
lower-dimensional tori under perturbation, that is for $\e\neq0$.
The Hamiltonian (\ref{eq:3.18}) can be interpreted as a
Hamiltonian of the form (\ref{eq:8.6}) with $\l_{i}=O(\e)$.
In that case one usually say that the lower-dimensional tori
are \textit{parabolic}. However, we can also interpret
the persisting tori, in the case $a_{i}>0$ and $\e<0$,
as \textit{degenerate hyperbolic tori}, where
``degenerate'' refers to the fact that the normal frequencies
vanish at $\e=0$.

Of course if the matrix $\partial_{\bbb}^{2}f_{\zerooo}(\bbb_{0})$
is negative definite, the same result of persistence of
hyperbolic invariant tori holds for $\e>0$. The case of
indefinite (i.e. neither positive nor negative defined)
matrices will be considered at the end of Section \ref{sec:10}.

\zerarcounters
\section{Summation of the divergent series --
dissipative systems}
\label{sec:9}

\noindent
The discussion of the model (\ref{eq:3.26}) introduced in
Section \ref{sec:3.4} proceeds very closely to the case
of the hyperbolic tori of Section \ref{sec:8}.
In that case $n=1$, hence both the propagators and the
self-energies are scalar. One has, formally,
\be
\MM^{[n]}(x,\e) = - (\rmi \e x)^{2} +
\tilde \MM^{[n]}(x,\e) , \qquad
\tilde \MM^{[n]}(0,\e) = a \e + O(\e^{2}) ,
\qquad a:= \partial_{x} g(c_{0}) \neq 0 . 
\label{eq:9.1} \ee
If for all $n'<n$ the dressed propagators
$G^{[n']}(\omt\cdot\nut,\e)$ can still be bounded proportionally
to $|x|^{-1}$ (as the undressed ones), the terms $O(\e^{2})$
are defined by a convergent series, so that
\be
\tilde \MM^{[n]}(x,\e) =
\tilde \MM^{[n]}(0,\e) + x \int_{0}^{1} {\rm d}t \,
\partial_{x} \tilde \MM^{[n]}(t x,\e) , \qquad
\partial_{x} \tilde \MM^{[n]}(x,\e) = O(\e^{2}) ,
\label{eq:9.2} \ee
and hence
\be
\d_{0}(x) - \MM^{[n-1]}(x,\e) =
\rmi x - \MM^{[n-1]}(x,\e) = \rmi x \left( 1 +  \rmi \e x \right)
- a \e - O(\e^{2}) - O(\e^{2}x) , \qquad x = \omt\cdot\nut ,
\label{eq:9.3} \ee
and hence
\be
\left| G^{[n]}(x,\e) \right| \le
\left| \d_{0}(x) - \MM^{[n-1]}(x,\e) \right| \le \frac{2}{|x|} ,
\label{eq:9.4} \ee
for $\e$ small enough. Actually one can prove that
$\tilde \MM^{[n]}(0,\e)$ is real for real $\e$ \cite{G5},
a property which becomes essential to deal with the case
in which the nondegeneracy condition $a \neq0$ is not
satisfied -- see Section \ref{sec:12.2} below.

Again the properties (\ref{eq:9.1}) and (\ref{eq:9.4}) are proved
together, by induction on $n$: for $n=0$ they trivially hold, and, by
assuming that they are satisfied up to $n-1$, one sees that the series
defining $\tilde \MM^{[n]}(x,\e)$ converge, and hence (\ref{eq:9.1})
can be proved for $n$; see \cite{GBD1,G4,G5} for details.

\begin{figure}[ht]
\vskip-.1truecm
\centering
\ins{286pt}{-58pt}{$\hbox{Re}\,\e$}
\ins{206pt}{-00pt}{$\hbox{Im}\,\e$}
\includegraphics[width=2in]{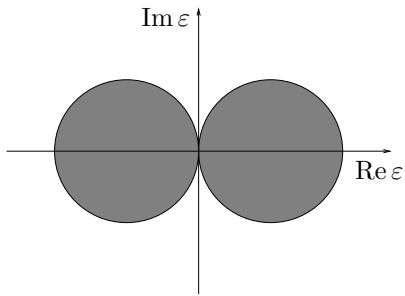}
\caption{Analyticity domain for the response solution.}
\label{fig:9.1}
\end{figure}

One can study the dependence of the response solution
on $\e$ in the complex domain. Of course one expects an
obstruction to analyticity along the imaginary axis
(so as it happened along the positive real axis
for the hyperbolic tori). In fact one can prove that
the response solutions are analytic in two disks
tangent to the imaginary axis at the origin -- see
Figure \ref{fig:9.1} --; of course only the disk to the right
is physically relevant, as it corresponds to $\e>0$
and hence to $\g>0$.

An interesting remark is that, while the periodic case
($d=p=1$) is trivial in the case of lower-dimensional tori,
this is no longer true for the model (\ref{eq:3.26}).
Indeed, if one takes $r=1$ in Section \ref{sec:3.3},
then no small divisors appear, so that the perturbation 
series is easily proved to be convergent (in particular
analyticity in $\e$ follows in that case). On the contrary
if one takes $d=m=1$ for the model (\ref{eq:3.26}),
one can still have arbitrarily large powers of $|\nut|$
because $\d_{1}(\omt\cdot\nut)/\d_{0}(\omt\cdot\nut)=\rmi\omt\cdot\nut$: then,
it is straightforward to see that one can construct trees
whose value grows like a factorial -- see Appendix \ref{app:C}.

This means that also in the case of periodic forcing,
the response solution to (\ref{eq:3.26}) is not analytic in $\e$.
However, one can prove that such a solution
is \textit{Borel summable} \cite{GBD1,GBD2};
we recall in Appendix \ref{app:D} the definition of
Borel summability -- see also \cite{Borel,Hardy,S}.
Note that an equation like (\ref{eq:3.26}) with periodic forcing
(and with $g(x)=x^{\mu}$ for $\mu\in[1.5,2.5]$)
naturally arises in electronic engineering,
and is known as the \textit{varactor equation} \cite{BDGM}.

It is proved in \cite{GBD2} that also in the case of
quasi-periodic forcing, the response solution turns out to
be Borel summable if $d=2$ and $\tau=1$ --  that is
in the case of frequency vectors with components whose ratio
is an irrational number of constant type \cite{Sch}.
A similar situation is encountered in the case of
hyperbolic tori: also in that case the function
$u(\ooo t,\e)$ is Borel summable if $d=2$ and $\tau=1$ \cite{CGGG}.

\zerarcounters
\section{Cantorisation -- elliptic tori}
\label{sec:10}

\noindent
Let us come back to the system of Section \ref{eq:3.3},
still assuming that the eigenvalues $a_{1},\ldots,a_{s}$ of
the matrix $\partial_{\bbb}^{2}f_{\zerooo}(\bbb_{0})$ are all positive,
but taking $\e>0$. In that cases, already for $n=-1$
one has -- see (\ref{eq:6.8}) and (\ref{eq:8.1}) --
\be
\MM^{[-1]}(x,\e) = \left( \begin{matrix} 0 & 0 \\
0 & \e\partial_{\bbb}^{2}f_{\zerooo}(\bbb_{0}) \end{matrix} \right) ,
\label{eq:10.1} \ee
so that the eigenvalues $x^{2} - \l^{[-1]}_{i}(x,\e)$ of the matrix
$\d_{0}(x)\,\id - \MM^{[-1]}(x,\e)$ are
\be
x^{2} - \l^{[-1]}_{i}(x,\e) = \begin{cases}
x^{2} , & i =1,\ldots, r , \\
x^{2} - \e a_{i-r} , & i=r+1,\ldots,d . \end{cases}
\label{eq:10.2} \ee
Hence, for fixed $\e$, we have problems for all $\nut\in\ZZZ^{r}$ such
that $\omt\cdot\nut$ is too close to some value $\pm\sqrt{a_{i}\e}$.
So, to give a meaning to $G^{[0]}(x,\e)$ we must require some
further Diophantine conditions, say
\be
\left| \left| \omt\cdot\nut \right| - \sqrt{\e a_{i}} \right| >
\g |\nut|^{-\tau'}
\quad \forall \nut\in\ZZZ^{r} \setminus\{\zerot\} ,
\label{eq:10.3} \ee
for some Diophantine exponent $\tau'\geq d$. This can be achieved
at the price of eliminating some values of $\e$.
Fixed $\e_{0}>0$ small enough, the subset $\gotE_{-1}'$
of values $\e\in[0,\e_{0}]$ for which all
Diophantine conditions (\ref{eq:10.3}) are satisfied has
large Lebesgue measure in $[0,\e_{0}]$, in the sense that
\be
\lim_{\e\to 0^{+}} \frac{{\rm meas}(\gotE_{-1}' \cap [0,\e])}{\e} = 1 ,
\label{eq:10.4} \ee
provided $\tau'$ is chosen larger enough than $\tau$,
say $\tau'>\tau+r$; see Appendix \ref{app:E}.
The property (\ref{eq:10.4}) can be stated
by saying that $\gotE_{-1}'$ has a Lebesgue density point at $\e=0$.

For all values $\e\in\gotE_{-1}'$ the propagators $G^{[0]}(x,\e)=
\Psi_{0}(x)\,(\d_{0}(x)\,\id - \MM^{[-1]}(x,\e))^{-1}$
can be formally defined. At this point, one could hope to
iterate the procedure. The main obstacle is that now
the dressed propagators are no longer bounded proportionally to
the undressed one: indeed it may happen that 
$x^{2} - \e a_{i-r}$ is much smaller than $x^{2}$. So we have
to modify the algorithm.

For simplicity, let us first reason once more by taking a
sharp decomposition as initially done in Section \ref{sec:5}.
Let us also assume, in the discussion below,
the self-energies to be well-defined:
we shall back later to such an issue.

We say that $\nut\ne 0$ is on scale $0$ if $|\omt\cdot\nut|\geq\g$ and
on scale $[\geq\!1]$ otherwise: for $\nut$ on scale $0$
we define $G^{[0]}(\omt\cdot\nut,\e)$ as in (\ref{eq:6.7}),
with $n=0$ and $\MM^{[-1]}(x,\e)$ given in (\ref{eq:10.1}). 
Given $\nu$ on scale $[\geq\!1]$ we say that $\nu$ is on scale $1$ if
$\min_{i=1,\ldots,d} |(\omt\cdot\nut)^{2}-
\l^{[-1]}_{i}(\omt\cdot\nut,\e)| \geq (2^{-1}\g)^{2} $,
and on scale $[\geq\!2]$ if $\min_{i=1,\ldots,d} |
(\omt\cdot\nut)^{2}-\l^{[-1]}_{i}(\omt\cdot\nut,\e) | < (2^{-1}\g)^{2}$.
For $\nu$ on scale $1$ we write $G^{[1]}(x,\e)$ as in (\ref{eq:6.7}),
with $n=1$ and $\MM^{[0]}(x,\e)$ written according to (\ref{eq:6.8}).
Call $\l^{[0]}_{i}(x,\e)$ the eigenvalues of $\MM^{[0]}(x,\e)$:
given $\nu$ on scale $[\geq\!2]$ we say that $\nu$ is on scale $2$ if
$\min_{i=1,\ldots,d} |(\omt\cdot\nut)^{2}-
\l^{[0]}_{i}(\omt\cdot\nut,\e)| \geq (2^{-2}\g)^{2} $,
and on scale $[\geq\!3]$ if $\min_{i=1,\ldots,d}|(\omt\cdot\nut)^{2}-
\l^{[0]}_{i}(\omt\cdot\nut,\e)|< (2^{-2}\g)^{2} $. And so on: eventually
we impose infinitely many Diophantine conditions, i.e.
\be
\left| \left| \omt\cdot\nut \right| -
\sqrt{ |\l^{[n]}_{i}(\omt\cdot\nut,\e)| } \right| >
2^{-(n+1)/2} \g |\nut|^{- \tau'} \quad
\forall \nut \in \ZZZ^{r} \setminus \{\zerot\} 
\label{eq:10.5} \ee
for all $i=1,\ldots,d$ and $n\ge -1$.

Even if we are successful in imposing the conditions (\ref{eq:10.5}),
the argument above is still incomplete, and it needs
a further modification. In order to bound the tree values
and the self-energies we need a bound on the number of lines
of fixed scale, in the spirit of the Siegel-Bryuno lemma.
This requires to compare the propagators of the lines
entering and exiting clusters which are not self-energy clusters
-- see the discussion at the end of Appendix \ref{app:F} -- and
this leads to further Diophantine conditions,
\be
\left| \left| \omt\cdot
\left( \nut_{1}\! - \!\nut_{2} \right) \right| \pm
\sqrt{ |\l^{[n]}_{i}(\omt\!\cdot\!\nut_{1},\e)| } \pm
\sqrt{ |\l^{[n]}_{j}(\omt\!\cdot\!\nut_{2},\e)| } \right| >
2^{-(n+1)/2} \g |\nut_{1} \!-\! \nut_{2}|^{- \tau'} \quad
\forall \nut_{1} \!\neq\! \nut_{2} \in
\ZZZ^{r} \!\setminus\! \{\zerot\} ,
\label{eq:10.6} \ee
for all $i,j=1,\ldots,d$ and $n\ge -1$.
For any fixed $\nut$ this would mean to impose the conditions
for all $\nut_{1}$ and $\nut_{2}$ such that $\nut_{1}-\nut_{2}=\nut$.
Unfortunately, these conditions are infinitely many, and for
all of them we would eliminate intervals of the same size:
therefore we would left with a zero measure set.
Of course, at the step $n=-1$, there would be no difficulty,
since the eigenvalues $\l^{[-1]}_{i}(x,\e)$ are independent
of $x$ -- see (\ref{eq:10.2}) --, but already at the step $n=0$
problems would arise.

So, instead of (\ref{eq:10.5}) and (\ref{eq:10.6}),
we can try to impose the Diophantine conditions
\begin{subequations}
\begin{align}
& \left| \left| \omt\cdot\nut \right| -
\sqrt{ \underline{\l}^{[n]}_{i}(\e) } \right| >
2^{-(n+1)/2} \g |\nut|^{- \tau'} \quad
\forall \nut \in \ZZZ^{r} \setminus \{\zerot\} ,
\label{eq:10.7a} \\
& \left| \left| \omt\cdot\nut \right| \pm
\sqrt{ \underline{\l}^{[n]}_{i}(\e) } \pm
\sqrt{ \underline{\l}^{[n]}_{j}(\e) } \right| >
2^{(n+1)/2} \g |\nut|^{- \tau'}
\quad \forall \nut \in \ZZZ^{r} \setminus \{\zerot\} ,
\label{eq:10.7b}
\end{align}
\label{eq:10.7} \end{subequations}
\vskip-.3truecm
\noindent for all $i,j=1,\ldots,d$ and $n\ge 0$, for suitable
numbers $\underline{\lambda}^{[n]}_{i}(\e)$ independent of $\nut$.
The advantage of (\ref{eq:10.7}) with respect to (\ref{eq:10.6})
is that for any $n$ we have to impose
that the quantities $\omt\cdot\nut$ are far enough only
from a finite number of values, that is $d$ values
for (\ref{eq:10.7a}) and $\le 4d^{2}$ values for (\ref{eq:10.7b}).

Thus, already for $n=-1$ we have to impose, besides the
conditions (\ref{eq:10.3}), also the conditions (\ref{eq:10.7b}).
This leaves a subset $\gotE_{-1}\subset \gotE_{-1}'$.
To prove that the set $\gotE_{-1}$ has still a Lebesgue density point
at $\e=0$, we need a lower bound on all the derivatives
${\rm d}(\sqrt{ \e a_{i}} \pm \sqrt{\e a_{j}})/{\rm d}\e$,
with $i\ne j$ when the sign minus is considered. This is easily
obtained if we assume that the eigenvalues $a_{1},\ldots,a_{s}$
are distinct, i.e. that there exists $a_{0}>0$ such that $|a_{i}-a_{j}|>
a_{0}$ for all $1 \le i < j \le s$; see Appendix \ref{app:E}.
Of course, this provides a further
assumption on the function $f_{\zerooo}(\bbb)$.

To deal with the cases $n\ge 0$ we define, iteratively,
\be
\underline{\l}^{[n]}_{i}(\e) = \begin{cases}
0 , & i =1,\ldots,r , \\
\lambda_{i}^{[n]}(\sqrt{\underline{\l}^{[n-1]}_{i}(\e)},\e) ,
& i = r+1 , \ldots,d . \end{cases}
\label{eq:10.8} \ee
In this way we obtain both that the eigenvalues $x^{2}-\l^{[n]}_{i}(x,\e)$
are bounded in terms of the quantities $x^{2}-\underline{\l}^{[n]}_{i}(\e)$
and that the sequences $\{\underline{\l}^{[n]}_{i}(\e)\}_{n=-1}^{\io}$
converge exponentially fast for all $i=r+1,\ldots,d$, that is
$|\underline{\l}^{[n]}_{i}(\e)-\underline{\l}^{[n-1]}_{i}(\e)|
\le K_{1} \rme^{-\ka_{1}2^{n/\tau'}}\e^{2}$ for suitable positive
constants $K_{1}$ and $\ka_{1}$; see Appendix \ref{app:F}.
Again, in order to impose the conditions (\ref{eq:10.7b})
we need a lower bound on the derivatives ${\rm d}(\sqrt{
\underline{\l}^{[n]}_{i}(\e) } \pm \sqrt{\underline{\l}^{[n]}_{j}(\e)
})/{\rm d}\e$, but these can be discussed as in the
case $n=-1$; see Appendix \ref{app:E}.

The discussion above is correct as far as the self-energies are
well-defined -- which we have simply assumed to be for the moment.
For instance, only if this is the case, when we write
$\underline{\l}_{i}^{[n]}(\e)=a_{i-r}\e +O(\e^{2})$
for $i=r+1,\ldots,d$, we can really say that the high order terms
are negligible with respect to the liner ones.
As in the case of maximal and hyperbolic tori,
we prove by induction that the self-energies are well defined.
To this aim, we need bounds on the number  of lines on scale $n$:
one can prove that bounds of the form (\ref{eq:5.8}) and (\ref{eq:7.2})
still holds, but with $\tau'$ instead of $\tau$  -- since the Diophantine
conditions involve the Diophantine exponent $\tau'$;
see Appendix \ref{app:F}. Up to this difference,
the strategy of the inductive proof is exactly
as in the previous cases.

The Diophantine conditions (\ref{eq:10.7a}) and (\ref{eq:10.7b})
are known as the \textit{first Melnikov conditions} and
\textit{second Melnikov conditions}, respectively.
Each condition shrinks further the set of allowed values of $\e$:
if $\gotE_{n-1}$ is the set of allowed values
found at the step $n-1$, then imposing the conditions
(\ref{eq:10.7}) leaves a subset $\gotE_{n} \subset \gotE_{n-1}$.
At each step the set of values which are removed has measure
proportional to a common value times $2^{-(n+1)/2}$: it was to obtain
this exponential factor that a factor $2^{-(n+1)/2}$ was introduced
in (\ref{eq:10.7}). Thus, eventually one is left with a set
$\gotE_{\io}$, which is still of large measure.
A closer inspection of $\gotE_{\io}$ reveals that $\gotE_{\io}$
is a Cantor set (that is a perfect, nowhere dense set).

To make the argument above really rigorous one should take a
smooth decomposition, such as that considered in Section \ref{sec:5}.
Moreover, it turns out to be convenient to use
functions $\underline{\l}^{[n]}_{i}(\e)$ which are differentiable
(in the sense of Whitney \cite{Poe1}) in $\e$, so that, instead of the
minimum of the eigenvalues, one should consider
a smooth version of it -- see Appendix \ref{app:F} for details
(see also \cite{G2,GG2,G3}).

The assumptions on $\partial_{\bbb}^{2}f_{\zerooo}(\bbb_{0})$
can be weakened by requiring that the eigenvalues
$a_{1},\ldots,a_{s}$ are such that $a_{i} \neq 0$ for all
$i$ and $a_{i}-a_{j} \neq 0$ for all $i\neq j$. In this case
\textit{lower-dimensional tori of mixed type} can be proved to exist
\cite{GG2,G3}.

\zerarcounters
\section{Stability and uniqueness}
\label{sec:11}

\noindent
The quasi-periodic solutions describing the maximal tori are linearly
stable \cite{AKN}. An interesting problem is that of uniqueness
of solutions. In other words, one can wonder whether there are
other quasi-periodic solutions with the same rotation vector
$\ooo$ as the solution studied in the previous sections.
Despite the apparent simplicity of the problem,
a proof of uniqueness has been given only recently \cite{FGS}.

The case of lower-dimensional tori is more difficult.
In principle there could be other quasi-periodic solutions with
the same rotation vector, which either do not admit
any perturbation expansion or admit a different expansion
or, even admitting the same expansion, are different.
For instance, the resummed expansion (\ref{eq:6.6}) \textit{a priori}
depends on the particular way the multiscale analysis is implemented,
and by slightly changing the procedure one could obtain a
different solution: this would imply infinitely many
solutions which have the same formal perturbation series.
In the case $d=2$ and $\tau=1$ all such functions coincide,
because they are all Borel summable -- see the last remark
in Appendix \ref{app:D} -- but in general there is no reason
why this should happen. At the present moment the problem
of uniqueness is still open.

In the case of the dissipative systems of Section \ref{sec:3.4}
one expects, on the ground of physical considerations,
the response solution to be either attractive or repulsive (which
means attractive for the time-reverted dynamics). More precisely,
under the further assumption that $\partial_{x} g(c_{0}) >0$
the response solution is expected to be asymptotically stable.
Indeed, this is what happens.

The proof -- very easy -- proceeds
as follows \cite{BDG1}. The analysis of Section \ref{sec:9}
shows that there exists a response solution $x_{0}(t)=x(\omv t,\e)=
c_{0}+O(\e)$. If we look for solutions of the form $x=x_{0}+\xi$:
then $x(t)\to x_{0}(t)$ as $t\to\io$
(i.e. $x_{0}$ is attracting) if and only
if $\xi(t)\to 0$ as $t\to\io$. The function $\xi$ must solves
the differential equation
\be
\ddot \xi + \g \dot \xi + P(\xi,x_{0}(t)) = 0 ,
\qquad P(\xi,x)=g(x+\xi)-g(x) = \partial_{x}g(x)\,\xi + O(\xi^{2})  ,
\label{eq:11.1} \ee
which can be rewritten as a system of first order equations
\be
\begin{cases} \dot \xi = y , & \\
\dot y = - \g y - P(\xi,x_{0}(t)) . \end{cases}
\label{eq:11.2} \ee
If we define $R(\xi,t)=P(\xi,x_{0}(t))/P(\xi,c_{0})$ we have
$R(0,t) = 1 + O(\e)$, so that $1/2 \le R(\xi,t) \le 2$
uniformly in $t$ and $\xi$, for $\e$ and $\xi$ small enough.
Then we can rescale time and variables by setting
\be
\tau(t) = \int_{0}^{t} {\rm d}t' \, \sqrt{R(\xi(t'),t')} , \qquad
\xi(t) = v(\tau(t)) , \qquad
y(t) = \sqrt{R(\xi(t),t)} w(\tau(t)) ,
\label{eq:11.3} \ee
which transforms the system (\ref{eq:11.2}) into
\be
\begin{cases}
v' = w , & \\
w' = - \g(v,t) \, w - P(v,c_{0}) , \end{cases}
\label{eq:11.4} \ee
with the prime denoting differentiation with respect to
time $\tau$ and
\be
\g(v,t) := \frac{1}{\sqrt{R}}
\left( \g + \frac{R'}{2\sqrt{R}} \right) .
\label{eq:11.5} \ee
If we neglect the friction term $\g(v,t)\,w$ in (\ref{eq:11.4}) we
obtain an autonomous system with constant of motion
\be
H(v,w) = \frac{1}{2} w^{2} + \int_{0}^{v} {\rm d}v'\, P(v',c_{0})  =
\frac{1}{2} w^{2} + \frac{1}{2} \partial_{x}g(c_{0}) \, v^{2} +
O(v^{3}) ,
\label{eq:11.6} \ee
so that the origin is a stable equilibrium point. Moreover
$\g(v,t) >0$ for $\e$ small enough (recall that $\g=1/\e$),
in a neighbourhood $U$ of the origin. Hence we can apply
Barbashin-Krasokvsky's theorem \cite{Kr} (or Lasalle's invariance principle
\cite{HK}) to conclude that the origin is
asymptotically stable and $U$ is contained in its basin of attraction.

In some cases, for instance if $g(x)=x^{2p+1}$, $p\in\NNN$,
and $f_{\zerooo} \neq 0$ (so that $\partial_{x}g(c_{0}) > 0$),
the response solution turns out to be a global attractor \cite{BDG1},
but of course in general it is only locally attracting.

The problem of uniqueness, mentioned about the lower-dimensional
tori, can be addressed also as to such response solutions.
Of course, the local attractiveness of the solution implies
local uniqueness. In other words, the response solution
$x(\omv t,\e)$ is the only quasi-periodic solution which
reduces to $c_{0}$ as $\e\to 0$.

\zerarcounters
\section{Generalisations}
\label{sec:12}

\noindent
In this last section, we review some possible directions one can follow
to generalise the results described in the previous sections.
Some of the these generalisations are discussed in the literature,
other are still open problems.

\subsection{Weaker Diophantine conditions}
\label{sec:12.1}

\noindent
Instead of the standard Diophantine condition (\ref{eq:1.5}) one can
consider weaker conditions, such as the \textit{Bryuno condition}
\cite{B}: a vector $\omt$ is said to satisfy
the Bryuno condition if $\calB(\omt) < \io$, where
\be
\calB(\omt) = \sum_{n=1}^{\io} \frac{1}{2^{n}}
\log  \frac{1}{\al_{n}(\omt)} , \qquad \al_{n}(\omt) = \inf
\{ |\omt\cdot\nut| : \nut \in \ZZZ^{d}
\hbox{ such that } 0<|\nut|\le 2^{n} \} .
\label{eq:12.1} \ee
All the results of the previous sections can be extended
to rotation vectors satisfying the Bryuno condition:
see \cite{G3} for maximal and lower-dimensional tori,
\cite{BeG1} for the standard map,
and \cite{G4} for dissipative systems.

For $d=2$ one can write $\omt=(\o_{1},\o_{2})=(1,\al)\o_{1}$,
where $\al=\o_{2}/\o_{1}$ is the rotation number. One can define
the \textit{Bryuno function}
\be
B(\al) = \sum_{n=1}^{\io} \frac{1}{q_{n}} \log q_{n+1} ,
\label{eq:12.2} \ee
where $\{q_{n}\}_{n\in\ZZZ}$ are the denominators of
the \textit{best approximants} of $\al$ \cite{Sch}.
Then the function $\calB(\omt)$ is equivalent
to the Bryuno function $B(\al)$, in the sense that one has
$C^{-1} < \calB(\omt)/B(\al) < C$ for a universal constant $C$ \cite{G3}.

An open problem is whether such a condition can be further weakened.
For $d=2$ the Bryuno condition is optimal: indeed, in the case
of the standard map Davie \cite{D} proved that if the rotation
number $\o\in\RRR$ does not satisfy the Bryuno condition
then there is no invariant curve with that rotation number $\o$.
Such a result can be even strengthened
by saying that the radius of convergence $\rho(\o)$ of the conjugating
function and the function $B(\o)$ are such that
\be
C_{1} \rme^{- 2 B(\o)} < \rho(\o) < C_{2} \rme^{- 2B(\o)} ,
\label{eq:12.3} \ee
for suitable universal constants $C_{1}$ and $C_{2}$;
see \cite{D,BeG1} for a proof of the last statement.
The proof of the lower bound in (\ref{eq:12.3}) relies on
deeper cancellations than those discussed
in Section \ref{sec:7}; see \cite{BeG1} for details.

Another Diophantine condition considered in the literature is the
so-called \textit{R\"ussmann condition} \cite{R1,R2,Poe2}, which has a
somewhat intricate definition if compared to (\ref{eq:12.1}).  For
$d=2$ such a condition is equivalent to the Bryuno condition \cite{R2}.

\subsection{Degenerate perturbations}
\label{sec:12.2}

\noindent
The dissipative systems introduced in Section \ref{sec:3.4}
have been considered under the \textit{nondegeneracy condition}
that $\partial_{x}g(c_{0})\neq0$ if $g(c_{0})=f_{\zerov}$. Such a
condition can be removed, and the existence of a response solution
can be proved under the only condition that there exists a
zero $c_{0}$ of odd order to the equation $g(x)-f_{\zerov}=0$;
see \cite{G4}, where it is also proved that there is no response solution
reducing to $c_{0}$ as $\e\to0$ if $c_{0}$ is a zero of even order.

Also in the case of lower-dimensional tori one can think of
relaxing the nondegeneracy condition that the matrix $\partial_{\bbb}^{2}
f_{\zerooo}(\bbb_{0})$ be nonsingular. In full generality,
this case is very hard. The case $s=1$ is already
nontrivial. In that case it has be proved that at least
one lower-dimensional torus always persists \cite{Che}.

A first difference with respect to the nondegenerate case
considered in Section \ref{sec:3.3} is that a formal power
series in $\e$ does not exist any longer. In \cite{GGG}
it is shown that, at least in some cases, a fractional
power series in $\e$ can be envisaged. The situation
is somewhat reminiscent of what happens
in \textit{Melnikov theory} \cite{Me1,CH,GH},
when the \textit{subharmonic Melnikov function} has degenerate zeroes
of odd order: in that case subharmonic solutions exist and are
analytic in a fractional power of $\e$ (\textit{Puiseux series}
\cite{Pui,BK}) -- cf. \cite{GBD3,CG}; see also \cite{ALGM,Pe}
for a similar situation in the case of limit cycles.
In the case of lower-dimensional tori,
the fractional power series in $\e$ do not converge,
but they can be resummed in order to give well-defined functions --
see \cite{GGG} for a complete discussion.

\subsection{More general systems}
\label{sec:12.3}

\noindent
One can also consider ordinary differential equations
more general than (\ref{eq:1.1}), say of the form
\be
D_{\e} u = F_{0}(u) + \e F(u, \omv t) ,
\label{eq:12.4} \ee
with $F_{0}$ real analytic. In that case one still assumes that
the unperturbed equation $D_{0}(u)=F_{0}(u)$ admits
a quasi-periodic solution $u_{0}(\ooo t)$.

The most general formulation of KAM theorem is within
this class -- see the Hamiltonian  (\ref{eq:3.2}).
The analysis performed in the previous sections for the
simplified Hamiltonians (\ref{eq:3.1}) and (\ref{eq:3.18})
can be extended to deal with these Hamiltonians;
we refer to \cite{GM2,G3} for details.
If the perturbation depends explicitly on the action variables,
then the formal solubility of the Hamilton equations relies on
a nondegeneracy condition of the unperturbed Hamiltonian, such as
the invertibility of the matrix $\partial_{\AAA}^{2}\calH_{0}(\AAA)$
(\textit{anisochronous condition}).
However, the KAM theorem can be extended also to
\textit{isochronous systems} with $\calH_{0}(\AAA) = \ooo\cdot\AAA$,
by assuming a nondegeneracy condition on the perturbation,
for instance that the matrix
\be
f_{\zerooo}(\AAA) = \int_{\TTT^{d}} \frac{{\rm d}\aaa}{(2\p)^{d}}
\, f(\aaa,\AAA) 
\label{eq:12.5} \ee
is invertible \cite{Ga2}. A proof along the lines of the
previous sections passes through the so-called
\textit{translated torus theorem} \cite{Mo2}
(also knows as \textit{theorem of the modifying terms} or
\textit{theorem of the counterterms}), which
says that, for any analytic function $f\!:\TTT^{n}\times\RRR^{n}$
and any Diophantine vector $\ooo\in\RRR^{n}$, there exists
a vector $\mmm(\e,\ooo)$ analytic in $\e$, such that the equations
\be
\begin{cases}
\dot\aaa = \ooo + \e \partial_{\AAA}f(\aaa,\AAA) +
\mmm(\e,\ooo) , & \\
\dot\AAA = - \e \partial_{\aaa}f(\aaa,\AAA) , & \end{cases}
\label{eq:12.6} \ee
admit a quasi-periodic solution with rotation vector $\ooo$
which is analytic in $\e$. A proof of such a theorem by using
the tree formalism can be found in \cite{BaG}. The theorem of
the translated torus
can be also formulated in Cartesian coordinates;
in that case the cancellation mechanisms leading to the
convergence of the series work in a rather different way \cite{CGP}.

In Section \ref{sec:3.1} we have considered only analytic
Hamiltonians. A more general formulation of the KAM theorem
requires only finite smoothness \cite{Mo1,Poe1}.
In certain cases, the tree formalism can be extended to nonanalytic
systems, such as some quasi-integrable systems of the form
(\ref{eq:3.3}) with $f$ in a class of $C^{p}$ functions
for some finite $p$ \cite{BGGM1,BGGM2}. However, up to exceptional cases,
the method described here seems to be intrinsically suited in cases
in which the vector fields are analytic -- from a physical point
of view this a quite reasonable assumption. The reason is that in
order to exploit the expansion (\ref{eq:2.3}), we need that $F$
be infinitely many times differentiable and we need a bound
on the derivatives. It is a remarkable property that, as shown
in Sections \ref{sec:8} and \ref{sec:9},  the
perturbation series can be given a meaning also in cases
where the solutions are not analytic in the perturbation parameter.

Equations of the form (\ref{eq:12.4}) also arise in problems
of electronic engineering and theory of circuits, usually
with periodic forcing. Such systems are resistive and hence
intrinsically dissipative. As examples one can consider
the \textit{saturating inductor circuit},
described by the equation
\be
 G(\dot x) \, \ddot x + \beta x + \e \g \dot x = \e f(\omv t) ,
\qquad G(v) = \frac{\al + v^{2}}{1+v^{2}} , \qquad
\al > 1, \quad \b>0 , \quad \g>0 ,
\label{eq:12.7} \ee
and the \textit{resonant injection-locked frequency divider},
described by the equation
\be
\begin{cases}
x' = \al y + \b x \left( 1 - x^{2} \right) +
\e x \left( 1 -x^{2} \right) f(\omv t) , & \qquad \al > \b > 1 , \\
y' = -x -y , & \end{cases}
\label{eq:12.8} \ee
with $m=1$ (i.e. periodic forcing) in both cases.
For both equations the dynamics at $\e=0$ is known:
the first system (\ref{eq:12.7}) admits a constant of motion
(although it is not Hamiltonian),
while the second one has a globally attracting limit cycle.
By following the same approach as described in Sections
\ref{sec:2} and \ref{sec:4}, one can study for $\e\neq0$
the existence of periodic solutions whose period is rational
with respect to the period of the forcing -- with
the major simplification with respect to the previous analysis
that no small divisors appear.
More precisely one can study the existence of subharmonic
solutions for the equation (\ref{eq:12.7}),
and the frequency locking phenomenon
for the equation (\ref{eq:12.8}). We refer to the literature
\cite{BDG2,BDG3} for details and results.

\subsection{Partial differential equations}
\label{sec:12.4}

\noindent
Finally, the analysis developed so far for ordinary differential
equations, can be extended to partial differential equations,
such as the \textit{nonlinear wave equation}
\be
\partial_{tt} u - \partial_{xx} u + \mu u = u^{3} ,
\qquad x \in [0,\p] , \qquad \mu \ge 0 ,
\label{eq:12.9} \ee
and the \textit{nonlinear Schr\"odinger equation}
\be
\rmi \partial_{t} u - \partial_{xx} u + \mu u = |u|^{2} u ,
\qquad x \in [0,\p] , \qquad \mu \ge 0 ,
\label{eq:12.10} \ee
with periodic or Dirichlet boundary conditions.
There exists a very wide literature about periodic and
quasi-periodic small amplitude solutions to nonlinear
one-dimensional partial differential equations such as
(\ref{eq:12.9}) and (\ref{eq:12.10}), starting from the
seminal work by Kuksin, Craig and Wayne \cite{Ku,W,CW}.
Recently results have been
obtained also in higher space dimension \cite{Bour1,Bour2,EK},
that is for $x\in [0,\p]^{D}$, $D>1$,
with periodic boundary conditions. 

By using the tree formalism, small amplitude periodic solutions
have been proved to exist, in dimension $D=1$,
both in the nonresonant case  -- $\mu$ in
a suitable Cantor set \cite{GM4} -- and
in the resonant case -- $\m=0$ \cite{GMP,GP1}.
Results have been obtained also in the higher
space dimensional case $D>1$ \cite{GP2}.
We refer to the original papers for a precise
formulation of the results and the proofs.

\appendix

\zerarcounters
\section{Proof of the Siegel-Bryuno lemma}
\label{app:A}

\begin{figure}[ht]
\centering
\ins{136pt}{-042pt}{$\theta=$}
\ins{210pt}{-068pt}{$T$}
\ins{170pt}{-050pt}{$\ell_{0}$}
\ins{246pt}{-024pt}{$\ell_{1}$}
\ins{256pt}{-038pt}{$\ell_{2}$}
\ins{230pt}{-067pt}{$\ell_{s}$}
\ins{284pt}{-001pt}{$\theta_{1}$}
\ins{303pt}{-028pt}{$\theta_{2}$}
\ins{286pt}{-077pt}{$\theta_{s}$}
\includegraphics[width=2.0in]{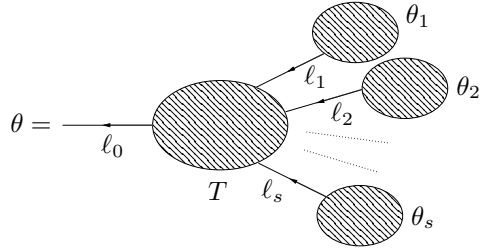}
\caption{Construction for the proof of the Siegel-Bryuno lemma.}
\label{fig:A.1}
\end{figure}

\noindent
The bound (\ref{eq:5.8}) follows from the fact that if $\gotN_{n}^{*}
(\theta)\neq0$ then $\gotN_{n}^{*}(\theta) \leq E(n,\theta):=
c\,K(\theta) 2^{-n/\tau}-1$, with $c=2^{2+1/\tau}$.
The last bound can be proved by induction on the order
$k(\theta)$ as follows. Given a tree $\theta$ let $\ell_{0}$ be its
root line, let $\ell_{1},\ldots,\ell_{s}$, $s\geq0$, be the lines on
scales $\ge n$ which are the closest to $\ell_{0}$, and let $\theta_{1},
\ldots, \theta_{s}$ the trees with root lines $\ell_{1},\ldots,\ell_{s}$,
respectively -- cf. Figure \ref{fig:A.1}. By construction all lines
$\ell$ in the subgraph $T$ have scales $n_{\ell}<n$, so that if
$n_{\ell_{0}}\ge n$ then $T$ is necessarily a cluster. Moreover,
all trees $\theta_{1},\ldots,\theta_{s}$ have orders strictly
less than $k(\theta)$, so that, by the inductive hypothesis,
for each $i=1,\ldots,s$ one has either $\gotN_{n}^{*}(\theta_{i})
\leq E(n,\theta_{i})$ or $\gotN_{n}^{*}(\theta_{i})=0$.

If either $\ell_{0}$ is not on scale $n$ or it
is on scale $n$ but exits a self-energy cluster then
$\gotN_{n}^{*}(\theta)=\gotN_{n}^{*}(\theta_{1})+\ldots+
\gotN_{n}^{*}(\theta_{s})$ and the bound $\gotN_{n}^{*}(\theta)
\leq E(n,\theta)$ follows by the inductive hypothesis. If $\ell_{0}$
does not exit a self-energy cluster and $n_{\ell_{0}}=n$
then $\gotN_{n}^{*}(\theta)=1+\gotN_{n}^{*}(\theta_{1})+\ldots+
\gotN_{n}^{*}(\theta_{s})$, and the lines $\ell_{1},\ldots,\ell_{s}$
enter a cluster $T$ with $K(T)=K(\theta)-(K(\theta_{1})+
\ldots+K(\theta_{s}))$. If $s\ge2$ the bound
$\gotN_{n}^{*}(\theta) \leq E(n,\theta)$ follows once more by the
inductive hypothesis. If $s=0$ then $\gotN_{n}^{*}(\theta)=1$;
on the other hand for $\ell_{0}$ to be on scale $n_{\ell_{0}}=n$ one
must have $|\omt\cdot\nut_{\ell_{0}}| < 2^{-n+1}\g$,
which, by the Diophantine condition (\ref{eq:1.8}), implies
$K(\theta) \ge |\nut_{\ell_{0}}|>2^{(n-1)/\tau}$, hence $E(n,\theta)>1$.
If $s=1$ call $\nut_{1}$ and $\nut_{2}$ the momenta of the lines
$\ell_{0}$ and $\ell_{1}$, respectively (in particular $\nu_{1}=
\nu_{\ell_{0}}$). By construction $T$ cannot
be a self-energy cluster, hence $\nut_{1}\neq\nut_{2}$.
Thus, by the Diophantine condition (\ref{eq:1.5}), one has
\be
2^{-n+2} \g \ge
\left| \omt\cdot\nut_{1} \right| + \left| \omt\cdot\nut_{2} \right| \geq
\left| \omt\cdot(\nut_{1}-\nut_{2}) \right| >
\frac{\g}{|\nut_{1}-\nut_{2}|^{\tau}} ,
\label{eq:A.1} \ee
because $n_{\ell_{0}}=n$ and $n_{\ell_{1}}\geq n$, and hence
\be
K(T) \geq \sum_{v\in N(T)} |\nut_{v}| \geq
|\nut_{1}-\nut_{2}| > 2^{(n-2)/\tau} ,
\label{eq:A.2} \ee
hence $T$ must contain ``many nodes''. In particular, one finds
also in this case $\gotN_{n}^{*}(\theta)=1+\gotN_{n}^{*}(\theta_{1})
\leq 1 + E(n,\theta_{1}) \le 1 + E(n,\theta) - c\,K(T)2^{-n/\tau}
\le E(n,\theta)$, where we have used that $c\,K(T)2^{-n/\tau} \ge 1$
by (\ref{eq:A.2}), provided $c=2^{2+1/\tau}$.

\begin{figure}[ht]
\centering
\ins{066pt}{-015.5pt}{$\theta=$}
\ins{098pt}{-022.5pt}{$\nut$}
\ins{150pt}{-060pt}{$T_{1}$}
\ins{193pt}{-060pt}{$T_{2}$}
\ins{296pt}{-060pt}{$T_{p}$}
\ins{142.5pt}{-022.5pt}{$\nut$}
\ins{188pt}{-022.5pt}{$\nut$}
\ins{244pt}{-022.5pt}{$\nut$}
\ins{294pt}{-022.5pt}{$\nut$}
\ins{336pt}{-044pt}{$\theta_{1}$}
\includegraphics[width=4.0in]{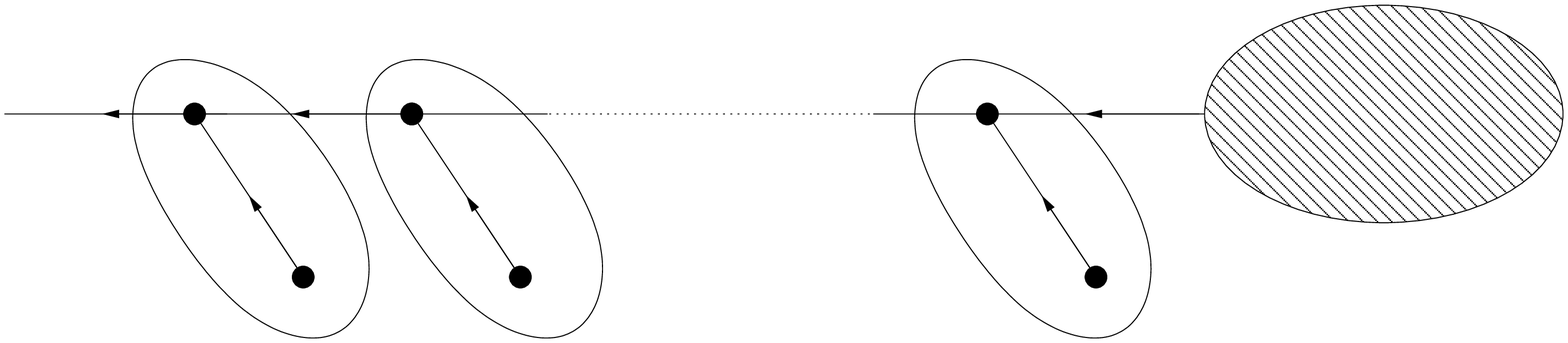}
\caption{Example of tree with accumulation of small divisors.}
\label{fig:A.2}
\end{figure}

The argument above shows that small divisors can accumulate only
by allowing self-energy clusters. That accumulation really occurs is
shown by the example in Figure \ref{fig:A.2}, where a tree $\theta$
of order $k$ containing a chain of $p$ self-energy clusters is depicted.
Assume for simplicity that $k/3$ is an integer: then if $p=k/3$
the subtree $\theta_{1}$ with root line $\ell$ is of order $k/3$.
If the line $\ell$ entering the rightmost self-energy cluster $T_{p}$
has momentum $\nut$, also the lines exiting the $p$ self-energy clusters
have the same momentum $\nut$. Suppose that $|\nut| \approx Nk/3$ and
$|\omt\cdot\nut|\approx \g/|\nut|^{\tau}$ (this is certainly possible
for some $\nut$). Then the value of the tree $\theta$ grows like
$a_{1}^{k}(k!)^{a_{2}}$, for some constants $a_{1}$ and $a_{2}$:
a bound of this kind prevents the convergence of the perturbation series.

\zerarcounters
\section{Siegel-Bryuno lemma for self-energies}
\label{app:B}

\noindent
We first prove (\ref{eq:7.3}). Call $\PP_{T}$ the path connecting
the entering and exiting lines of $T$. If $T\in \RR_{n}$
then $T$ contains at least a line on scale $\ge n$.
If there is one line $\ell\notin\PP_{T}$ on scale $\ge n$, one has
$K(T) \ge |\nut_{\ell}|$ and $\g|\nut_{\ell}|^{-\tau}< |\omt\cdot
\nut_{\ell}| \leq 2^{-n+1}\g$, so that $K(T) > 2^{(n-1)/\tau}$.
Otherwise, let $\ell$ be the line $\ell\in\PP_{T}$ on scale $\ge n$
which is closest to the line $\ell_{1}$ entering $T$.
Call $\tilde T$ the subgraph of $T$ consisting of all lines and
nodes of $T$ preceding $\ell$. By construction,
all lines in $\tilde T$ have scale $< n$, and hence
$\nut_{\ell} \neq \nut_{\ell_{1}}$, otherwise $\tilde T$
would be a self-energy cluster. Therefore one has
$\g|\nut_{\ell}-\nut_{\ell_{1}}|^{-\tau} < |\omt\cdot(\nut_{\ell}-
\nut_{\ell_{1}})| < 2^{-n+2}\g$, which yields $K(T) \ge
|\nut_{\ell}-\nut_{\ell_{1}}| > 2^{(n-2)/\tau}$.

To prove (\ref{eq:7.2}) one considers a more general class of graphs.
We say that $T$ is a graph of class $\SS(n,n')$ if it has
one entering line and one exiting line both on scale $\ge n'$
and all the lines contained in $T$ are on scale $\le n$.
Define $K(T)$ and $\gotN_{n}(T)$ as done in Section \ref{sec:7}
for the self-energy clusters. We want to prove the bound
$\gotN_{n'}^{*}(T) \leq \max\{2 K(T) 2^{(3-n)/\tau}-1,0\}$.
Let $\PP_{T}$ the path connecting the entering and exiting lines
of $T$, and let $N(\PP_{T})$ the set of nodes connected by
lines of $\PP_{T}$. If all the lines along $\PP_{T}$ have scale
$<n'$, then $\gotN_{n'}(T) = \gotN_{n'}(\theta_{1})+\ldots+
\gotN_{n'}(\theta_{m})$, where $\theta_{1},\ldots,\theta_{m}$
are the trees contained in $T$ which have the root in a
node $v\in N(\PP_{T})$. In that case the bound follows
from (\ref{eq:5.8}). If there exists a line $\ell\in \PP_{T}$
on scale $\ge n'$, call $T_{1}$ and $T_{2}$ the subgraphs of $T$
such that $L(T)=\{\ell\} \cup L(T_{1}) \cup L(T_{2})$.
Both $T_{1}$ and $T_{2}$ are of class $\SS(n,n')$, so that,
in the case in which both $T_{1}$ and $T_{2}$ contain lines 
on scale $\ge n'$, by the inductive hypothesis one finds
$\gotN_{n'}(T) \le 1 + \gotN_{n'}(T_{1}) + \gotN_{n'}(T_{2}) \le
2 K(T) 2^{(n-3)/\tau} - 1$. If $T_{1}$ contains no line
on scale $\ge n'$ then one realises that one must have
$K(T_{1}) > 2^{(n-2)/\tau}$, and the same holds for $T_{2}$,
so that the bound follows also in these cases.
Finally, (\ref{eq:7.2}) follows from the previous bound by
noting that a self-energy cluster $T\in\RR_{n}$ is a graph
of class $\SS(n,n')$ for all $n<n'$. 

\zerarcounters
\section{Accumulation of small divisors for dissipative systems}
\label{app:C}

\noindent
We want to construct for the model of Section \ref{sec:3.4}
a tree $\theta$ whose value $\Val(\theta)$ grows like a factorial.
Let $\theta$ be the tree with $k$ nodes $v_{1},\ldots,v_{k}$,
such that $s_{v_{i}}=k_{v_{i}}=1$ and $\rho_{v_{i}}=0$ for all
$i=1,\ldots,k-1$, while $s_{v_{k}}=0$ and $\rho_{v_{k}}=k_{v_{k}}=1$.
Let $\nut_{v_{k}}=\nut$ such that $|\omt\cdot\nut| \approx
\g/|\nut|^{\tau}$. The value of the labels $\rho_{v_{i}}$
for $i=1,\ldots,k-1$ implies that $\nut_{v_{i}}=0$ for $i=1,\ldots,k-1$,
and hence all the lines have the same momentum $\nut$.
Then one has $\Val(\theta)=
(\rmi\omt\cdot\nut)^{2(k-1)}f_{\nut}(\rmi\omt\cdot\nut)^{-k}=
(\rmi\omt\cdot\nut)^{k-2}f_{\nut}$, which can be bounded by
$(k-2)!(2/\xi)^{k}\Xi_{0}\rme^{-\xi|\nut|/2}$ for large $k$.

\zerarcounters
\section{Borel summability}
\label{app:D}

\noindent
Let $f(\e)=\sum_{n=1}^{\io} a_{n}\e^{n}$ a formal power series (which
means that the sequence $\{a_{n}\}_{n=1}^{\infty}$ is well-defined).
We say that $f(\e)$ is \textit{Borel summable} if
\begin{enumerate}
\item $B(p):=\sum_{n=1}^{\infty} a_{n} p^{n}/n!$
converges in some circle $|p|<\delta$,
\item $B(p)$ has an analytic continuation to a neighbourhood
of the positive real axis, and
\item $g(\e)=\int_{0}^{\io} {\rm e}^{-p/\e} B(p) \, {\rm d}p$
converges for some $\e>0$.
\end{enumerate}
Then the function $B(p)$ is called the \textit{Borel transform}
of $f(\e)$, and $g(\e)$ is the \textit{Borel sum} of $f(\e)$.
Moreover if the integral defining $g(\e)$ converges for
some $\e_{0}>0$ then it converges in the circle
$\hbox{Re}\,\e^{-1} > \hbox{Re}\,\e_{0}^{-1}$.

A function which admits the formal power series expansion $f(\e)$
is called Borel summable if $f(\e)$ is Borel summable;
in that case the function equals the Borel sum $g(\e)$.

A remarkable property of Borel summable functions is that
if two functions $f(\e)$ and $g(\e)$ are both
Borel summable and admit the same power series expansion,
then the two functions coincide.
 
\zerarcounters
\section{Excluded values of the perturbation parameter}
\label{app:E}

\noindent
Set $a=\min\{a_{1},\ldots,a_{s}\}$ and $A=\max\{a_{1},\ldots,a_{s}\}$.
In order to impose the Diophantine conditions (\ref{eq:10.3})
we have to exclude all values of $\e\in[0,\e_{0}]$ such that
$\left| |\omt\cdot\nut| - \sqrt{\e a_{i}}\right|\le \g|\nut|^{-\tau'} $
for some $i=1,\ldots,s$ and some $\nut\neq\zerot$.
Of course, we can confine ourselves
to the vectors $\nut\in\ZZZ^{r}$ such that
$|\nut| \ge m_{0}:=(\g/4\sqrt{\e_{0}A})^{-1/\tau}$, because
one has $|\omt\cdot\nut|>4\sqrt{\e_{0}A}$ and hence
$\left| |\omt\cdot\nut| - \sqrt{\e a_{i}}\right| > 
|\omt\cdot\nut|/2 > \g|\nut|^{-\tau}/2$ for $|\nut| < m_{0}$.
For all $|\nut| \ge m_{0}$ we can introduce an interpolation parameter
$t\in[-1,1]$ by setting $|\omt\cdot\nut| - \sqrt{\e(t,\nut)\,a_{i}} =
t\,\g |\nut|^{-\tau'}$, so that
\be
\left| \frac{{\rm d}}{{\rm d}t} \e(t,\nut) \right| \le
\frac{\g}{|\nut|^{\tau'}} \frac{2\sqrt{\e(t,\nut)}}{\sqrt{a_{i}}} \le
\frac{\g}{|\nut|^{\tau'}} \frac{2\sqrt{\e_{0}}}{\sqrt{a}} 
\label{eq:E.1} \ee
for all $\e(t,\nut)\in\gotE_{-1}'$.
Therefore we have to exclude a set $\gotE_{-1}' \subset [0,\e_{0}]$
of measure
\be
{\rm meas}(\gotE_{-1}') = \int_{\gotE_{-1}'} \!\!\!\! {\rm d}\e \le
\sum_{\substack{ \nut \in \ZZZ^{r} \\|\nut| \ge m_{0}}}
\!\!\!\! \int_{-1}^{1} \!\!\!\! {\rm d}t
\left| \frac{{\rm d}}{{\rm d}t} \e(t,\nut) \right| \le
\sum_{\substack{ \nut \in \ZZZ^{r} \\|\nut| \ge m_{0}}}
\frac{\g}{|\nut|^{\tau'}} \frac{4\sqrt{\e_{0}}}{\sqrt{a}} \le
C \g \left( \frac{\sqrt{\e_{0}A}}{\g} \right)^{(\tau'-r)/\tau}
\!\!\!\frac{\sqrt{\e_{0}}}{\sqrt{a}} ,
\label{eq:E.2} \ee
for some universal constant $C$. Hence
${\rm meas}(\gotE_{-1}')$ is much smaller than
$\e_{0}$ if $\e_{0}$ is small and $\tau'>\tau+r$.

To impose the Diophantine conditions (\ref{eq:10.7a}) one can
reason in the same way. One uses that $\underline{\l}^{[n]}_{i}(\e) =
a_{i-r}\e+O(\e^{2})$ for $i=r+1,\ldots,d$ (see Appendix \ref{app:F}),
which yields $|\underline{\l}^{[n]}_{i}(\e)| \le 2a_{i-r}\e$ and
$|{\rm d} \underline{\l}^{[n]}_{i}(\e)/{\rm d}\e| \ge a_{i-r}/2$; here
and henceforth the derivative is in the sense of Whitney.
Again we have to consider only the values $\nut\in\ZZZ^{r}$ such that
$|\nut|\ge m_{0}$. We define $\gotE_{n}'$ as the set
of values $\e\in[0,\e_{0}]$ which do not satisfy (\ref{eq:10.7a})
and for all $|\nut|\ge m_{0}$ and $t\in[-1,1]$ we write
$|\omt\cdot\nut| - \sqrt{\underline{\l}^{[n]}_{i}(\e(t,\nut))}
= t\,2^{-(n+1)/2}\g |\nut|^{-\tau'}$. Then
\be
\left| \frac{{\rm d}}{{\rm d}t} \e(t,\nut) \right| \le
2^{-(n+1)/2}
\frac{\g}{|\nut|^{\tau'}} \frac{4\sqrt{2\e}}{\sqrt{a_{i-r}}} \le
2^{-(n+1)/2} \frac{\g}{|\nut|^{\tau'}} \frac{4\sqrt{2\e_{0}}}{\sqrt{a}} 
\label{eq:E.3} \ee
for all $\e(t,\nut)\in\gotE_{n}'$, and hence
\be
{\rm meas}(\gotE_{n}') \le C 2^{-(n+1)/2} \g
\left( \frac{\sqrt{\e_{0}A}}{\g} \right)^{(\tau'-r)/\tau}
\frac{\sqrt{\e_{0}}}{\sqrt{a}} ,
\label{eq:E.4} \ee
for some universal constant $C$. Therefore one has ${\rm meas}
(\gotE_{n}') = 2^{-(n+1)/2} o(\e_{0})$, provided $\tau'>\tau+r$.

To impose the Diophantine conditions (\ref{eq:10.7b}), the only
difference is that we need a lower bound on the derivatives 
${\rm d} (\sqrt{\underline{\l}^{[n]}_{i}(\e)}\pm
\sqrt{\underline{\l}^{[n]}_{j}(\e)})/{\rm d}\e$,
$r+1\le i,j \le d$, with
$i\ne j$ when the sign minus is taken.
One easily realises that the conditions with the sign plus
do not present any further difficulty with respect to
the first Melnikov conditions. Moreover, for $\e$ small enough and
$r+1 \le i\ne j \le d$ one has
\be
\left| \frac{{\rm d}}{{\rm d}\e} \left(
\sqrt{\underline{\l}^{[n]}_{i}(\e)} -
\sqrt{\underline{\l}^{[n]}_{j}(\e)} \right) \right| \ge
\frac{a_{0} }{8\sqrt{A\e}} ,
\qquad a_{0} := \min_{1 \le k\ne h \le s} \left| a_{k}-a_{h} \right| .
\label{eq:E.5} \ee
To deduce (\ref{eq:E.5}) one uses that for $i=r+1,\ldots,d$.
one has $\l^{[-1]}_{i}(\e)=a_{i-r} \e$
and $\l^{[n]}_{i}(\e)=a_{i-r} \e +O(\e^{2})$, $n \ge 0$.
By defining $\gotE_{n}''$ as the set
of values $\e\in[0,\e_{0}]$ which do not satisfy (\ref{eq:10.7b}),
we obtain
\be
{\rm meas}(\gotE_{n}'') \le C 2^{-(n+1)/2} \g
\left( \frac{\sqrt{\e_{0}A}}{\g} \right)^{(\tau'-r)/\tau}
\frac{\sqrt{A\e_{0}}}{a_{0}} ,
\label{eq:E.6} \ee
for some universal constant $C$, so that once more
${\rm meas}(\gotE_{n}'')=2^{-(n+1)/2} o(\e_{0})$ for $\tau'>\tau+r$.
The sets $\gotE_{n}$ in Section \ref{sec:10} are defined
as $\gotE_{n}=\gotE_{n}'\cup \gotE_{n}''$.

\zerarcounters
\section{Multiscale analysis for elliptic tori}
\label{app:F}

\noindent
To extend the multiscale analysis to the case of elliptic tori,
we slightly change the recursive definition of propagators
and self-energies. We set $\Delta^{[-1]}(x,\e)=x^{2}$ and
\be
\Delta^{[n]}(x,\e) =
\left( \frac{1}{d} \sum_{i=1}^{d}
\left( x^{2} - \underline{\l}^{[n-1]}_{i}(\e) \right)^{-2}
\right)^{-1/2} \!\!\!\!\!\!\!\!\!\! , \qquad n \ge 0 , 
\label{eq:F.1} \ee
and define
\be
\Xi_{n}(x,\e) = \prod_{p=0}^{n} \chi_{p}(\Delta^{[p-1]}(x,\e)) ,
\qquad \Psi_{n}(x,\e) = \psi_{n}(\Delta^{[n-1]}(x,\e))
\prod_{p=0}^{n-1} \chi_{p}(\Delta^{[p-1]}(x,\e)) ,
\label{eq:F.2} \ee
with the functions $\chi_{n}$ and $\psi_{n}$ defined as in
Section \ref{sec:5}, with the only difference that in (\ref{eq:5.2})
$\g$ and $\g/2$ are replaced with $\g^{2}$ and $\g^{2}/4$, respectively.

In terms of the quantities (\ref{eq:F.1}) and (\ref{eq:F.2}),
the propagators $\GGG_{\ell}=G^{[n_{\ell}]}(\omt\cdot\nut_{\ell},\e)$
are defined iteratively as
\begin{subequations}
\begin{align}
& G^{[n]}(x,\e) =
\Psi_{n}(x,\e) \left( \d_{0}(x) \id -
\MM^{[n-1]}(x,\e) \right)^{-1}
\label{eq:A6.3a} \\
& \MM^{[n]}(x,\e) = \MM^{[n-1]}(x,\e) +
\Xi_{n}(x,\e) \, M^{[n]}(x,\e) , \qquad
M^{[n]}(x,\e) = \!\! \sum_{T\in\RR_{n}} \e^{k(T)} \, \VV_{T}(x) .
\label{eq:F.3b}
\end{align}
\label{eq:F.3} \end{subequations}
\vskip-.3truecm
\noindent Finally the numbers $\underline{\l}^{[n]}_{i}(\e)$
are defined according to (\ref{eq:10.8}), where $\l^{[n]}_{i}(x,\e)$
are the eigenvalues of the matrix $\MM^{[n]}(x,\e)$.

If a line $\ell$ is on scale $n$ and $\GGG_{\ell} \ne 0$, then one has
\begin{subequations}
\begin{align}
& \min_{1 \le i \le d} \left| (\omt\cdot\nut)^{2} -
\underline{\l}^{[p]}_{i}(\e) \right| \le 2^{-2 p} \g^{2} ,
\qquad 0 \le p \le n - 2 ,
\label{eq:F.4a} \\
& \min_{1 \le i \le d} \left| (\omt\cdot\nut)^{2} -
\underline{\l}^{[n-1]}_{i}(\e) \right| \ge \frac{1}{4\sqrt{d}}
2^{-2 n} \g^{2} .
\label{eq:F.4b}
\end{align}
\label{eq:F.4}
\end{subequations}
\vskip-.3truecm
\noindent Therefore, setting $\omt\cdot\nut=x$, if $x>0$ one has
\be
\left| \l^{[n-1]}_{i}(x,\e) - \underline{\l}^{[n-1]}_{i}(\e) \right|
\le \max_{x} \left| \partial_{x} \l^{[n-1]}_{i}(x,\e) \right| \,
\left| \sqrt{\underline{\l}^{[n-2]}_{i}(\e)} - x \right| ,
\label{eq:F.5} \ee
where $\partial_{x} \l^{[n-1]}_{i}(x,\e)=O(\e^{2})$ and
\be
\left| \sqrt{\underline{\l}^{[n-2]}_{i}(\e)} - x \right| \le
\frac{ | \underline{\l}^{[n-2]}_{i}(\e) - x^{2} |}
{| \sqrt{\underline{\l}^{[n-2]}_{i}(\e)} + x | } \le
\frac{2^{-2(n-2)}\g^{2}}{\sqrt{\e a}} ,
\label{eq:F.6} \ee
and hence
\be
\left| \l^{[n-1]}_{i}(x,\e) - \underline{\l}^{[n-1]}_{i}(\e) \right|
\le C \, \e_{0} 2^{-2n} ,
\label{eq:F.7} \ee
for some positive constant $C$. Therefore (\ref{eq:F.4b}) and
(\ref{eq:F.7}) imply
\be
\left| x^{2} - \l^{[n-1]}_{i}(x,\e) \right| \ge
\left| x^{2} - \underline{\l}^{[n-1]}_{i}(\e) \right| -
\left| \l^{[n-1]}_{i}(x,\e) - \underline{\l}^{[n-1]}_{i}(\e) \right|
\ge \frac{1}{2} \left| x^{2} - \underline{\l}^{[n-1]}_{i}(\e) \right| .
\label{eq:F.8} \ee
The case $x<0$ is discussed in the same way noting that
$\l^{[n-1]}_{i}(-x,\e)=\l^{[n-1]}_{i}(x,\e)$, because of (\ref{eq:7.1}).
Therefore the eigenvalues $x^{2} - \l^{[n]}_{i}(x,\e)$
can be bounded from below by half the quantities
$x^{2} - \underline{\l}^{[n]}_{i}(\e)$.

The property $| \underline{\l}^{[n]}_{i}(\e)- \underline{\l}^{[n-1]}_{i}
(\e)| \le C_{1} \rme^{-\ka_{1}2^{n/\tau'}} \e^{2}$
mentioned after (\ref{eq:10.8}) follows from
the expression (\ref{eq:F.3b}) for $\MM^{[n]}(x,\e)-\MM^{[n-1]}(x,\e)$,
the bound (\ref{eq:7.5}) for the values of the self-energy clusters,
and the bound $\sum_{v\in N(T)}|\nut_{v}| > c'' 2^{n/\tau'}$
which holds for any $T\in\RR_{n}$.
 
Finally we want to show that the bounds
\be
\gotN_{n}^{*}(\theta) \leq c\,2^{-n/\tau'} K(\theta) , \qquad
\qquad \gotN_{n'}^{*}(T) \leq c\,2^{-n'/\tau'} K(T) , \quad
T \in \RR_{n} , \quad n'\le n ,
\label{eq:F.9} \ee
hold with the multiscale analysis described above.
One proceed as in Appendices \ref{app:A} and \ref{app:B},
with the following changes. If the propagator of
a line with momentum $\nut$ and scale $n$ is non-zero,
then (\ref{eq:F.4}) imply $|\nut| \ge 2^{(n-2)/\tau'}$.
When discussing the analogous of the case $s=1$
in Appendix \ref{app:A}, then (\ref{eq:A.1}) must be replaced with
\begin{eqnarray}
2^{-n+3} \g & \!\!\! \ge \!\!\! &
| \, \omt\cdot\nut_{1} + \s_{1}
\sqrt{\underline{\l}^{[n-2]}_{i}(\e)} \, |
+ | \, \omt\cdot\nut_{2} + \s_{2}
\sqrt{\underline{\l}^{[n-2]}_{j}(\e)} \, |
\nonumber \\
& \!\!\! \ge \!\!\! &
| \, \omt\cdot(\nut_{1}-\nut_{2}) +
\s_{1} \sqrt{\underline{\l}^{[n-2]}_{i}(\e)} - 
\s_{2} \sqrt{\underline{\l}^{[n-2]}_{j}(\e)} \, | \ge
\frac{\g}{|\nut_{1}-\nut_{2}|^{\tau'}} ,
\label{eq:F.10} \end{eqnarray}
where the signs $\s_{1},\s_{2} \in \{\pm\}$ and the labels $i,j\in
\{1,\ldots,d\}$ are such that the first inequality is satisfied.
Analogously one discussed the case of the self-energy clusters.

In particular (\ref{eq:F.10}) explains why the second
Melnikov conditions are necessary. Of course, if we had used
directly the eigenvalues $\l^{[n]}_{i}(x,\e)$, instead of
the quantities $\underline{\l}^{[n]}_{i}(\e)$, we would
have require (\ref{eq:10.6}) instead of (\ref{eq:10.7b}).
We have already seen in Section \ref{sec:10} why this
was not allowed.

\vspace{.5truecm}
\noindent \textbf{Acknowledgments.} I thank L. Corsi
for her very careful comments on the manuscript.


\end{document}